\documentclass[11pt]{article}
\usepackage{epsfig,epsf,fancybox}
\usepackage{amsmath}
\usepackage{mathrsfs}
\usepackage{amssymb}
\usepackage{graphicx}
\usepackage{color}
\usepackage{multirow}
\usepackage{paralist}
\usepackage{verbatim}
\usepackage{galois}
\usepackage{algorithm}
\usepackage[toc,page,header]{appendix}
\usepackage{titletoc}

\floatstyle{ruled}
\newfloat{algorithm}{tbp}{loa}
\providecommand{\algorithmname}{Algorithm}
\floatname{algorithm}{\protect\algorithmname}
\usepackage{boxedminipage}
\usepackage{booktabs}
\usepackage{accents}
\usepackage{stmaryrd}
\usepackage[table]{xcolor}
\usepackage{hhline}
\DeclareSymbolFont{extraup}{U}{zavm}{m}{n}
\DeclareMathSymbol{\varheart}{\mathalpha}{extraup}{86}
\usepackage{natbib}
\usepackage{url}
\usepackage{tikz}
\usepackage{algorithmic}
\usepackage{xcolor}
\usepackage{microtype}
\usepackage{graphicx}
\usepackage{subfigure}
\usepackage{tcolorbox}
\usepackage{enumitem}
\usepackage{booktabs} 
\usepackage{multirow}
\usepackage{changepage}

\usepackage[colorlinks=true,breaklinks=true,bookmarks=true,urlcolor=teal,citecolor=teal,linkcolor=teal,bookmarksopen=false,draft=false]{hyperref}  
\usepackage[nameinlink, capitalize,noabbrev]{cleveref}
\crefname{assumption}{Assumption}{Assumptions}
\crefname{Assumption}{Assumption}{Assumptions}
\crefname{Lemma}{Lemma}{Lemmata}
\crefname{lemma}{Lemma}{Lemmata}
\crefname{Theorem}{Theorem}{Theorems}
\crefname{theorem}{Theorem}{Theorems}
\crefname{corollary}{Corollary}{Corollaries}
\crefname{proposition}{Proposition}{Propositions}
\crefname{claim}{Claim}{Claims}
\crefname{procedure}{Procedure}{Procedures}
\crefname{algorithm}{Algorithm}{Algorithms}
\crefname{figure}{Figure}{Figures}
\crefname{remark}{Remark}{Remarks}
\crefname{section}{Section}{Sections}
\crefname{procedure}{Procedure}{Procedures}
\crefname{definition}{Definition}{Definitions}
\crefname{example}{Example}{Examples}
\crefname{table}{Table}{Tables}
\crefname{equation}{}{}
\crefname{enumi}{}{}
\crefname{item}{Item}{Items}

\usepackage{ifthen}
\newboolean{Arxiv}
\setboolean{Arxiv}{true}   

\usepackage{geometry}
\numberwithin{equation}{section}

\title{\bf{\LARGE{Accelerating Low-Rank Factorization-Based Semidefinite Programming Algorithms on GPU}}}

\author{
Qiushi Han\thanks{University of Illinois Urbana-Champaign, \{\href{mailto:joshhan2@illinois.edu}{joshhan2@illinois.edu}, \href{mailto:hanwenl4@illinois.edu}{hanwenl4@illinois.edu}\}} \quad\quad\quad
Zhenwei Lin\thanks{Shanghai University of Finance and Economics, \href{mailto:zhenweilin@163.sufe.edu.cn}{zhenweilin@163.sufe.edu.cn}}  \quad\quad\quad
Hanwen Liu\footnotemark[1] \\
\\
Caihua Chen\thanks{Nanjing University, \href{mailto:chchen@nju.edu.cn}{chchen@nju.edu.cn}} \quad\quad\quad
Qi Deng\thanks{Shanghai Jiao Tong University, \{
\href{mailto:qdeng24@sjtu.edu.cn}{qdeng24@sjtu.edu.cn}, \href{mailto:ddge@sjtu.edu.cn}{ddge@sjtu.edu.cn}\}} \quad\quad\quad
Dongdong Ge\footnotemark[4] \quad\quad\quad
Yinyu Ye\thanks{Stanford University, \href{mailto:yinyu-ye@stanford.edu}{yinyu-ye@stanford.edu}} \quad
}

\usepackage{bm}

\global\long\def\inprod#1#2{\big\langle #1,#2\big\rangle }%

\global\long\def\inner#1#2{\langle#1,#2\rangle}%

\global\long\def\norm#1{\Vert#1\Vert}%
\global\long\def\Bnorm#1{\Big\Vert#1\Big\Vert}%

\global\long\def\brbra#1{\big(#1\big)}%
\global\long\def\Brbra#1{\Big(#1\Big)}%
\global\long\def\rbra#1{(#1)}%

\global\long\def\bsbra#1{\big[#1\big]}%
\global\long\def\abs#1{\vert#1\vert}%
\global\long\def\bcbra#1{\big\{#1\big\}}%
\global\long\def\vertiii#1{\left\vert \kern-0.25ex  \left\vert \kern-0.25ex  \left\vert #1\right\vert \kern-0.25ex  \right\vert \kern-0.25ex  \right\vert }%
\global\long\def\mcal#1{\mathcal{#1}}%
\global\long\def\mbb#1{\mathbb{#1}}%
\global\long\def\argmin{\operatornamewithlimits{argmin}}%
\global\long\def\nnz{\operatornamewithlimits{nnz}}%
\global\long\def\st{\operatornamewithlimits{s.t.}}%
\global\long\def\and{\mathrm{and}}%
\global\long\def\vep{\varepsilon}%
\global\long\def\Abf{\mathbf{A}}%
\global\long\def\vec{\operatorname{vec}}%
\global\long\def\mat{\operatorname{mat}}%

\newcommand{\reopt}{\texttt{re-opt }}
\date{July 2024}

\begin{document}
\maketitle

\begin{abstract}
In this paper, we address a long-standing challenge: \emph{how to achieve both efficiency and scalability in solving semidefinite programming problems.} We propose breakthrough acceleration techniques for a wide range of low-rank factorization-based first-order methods using GPUs, making the computation much more efficient and scalable. To illustrate the idea and effectiveness of our approach, we use the low-rank factorization-based SDP solver, LoRADS, as an example, which involves both the classic Burer-Monterio method and a novel splitting scheme with a starting logarithmic rank. Our numerical results demonstrate that the accelerated GPU version of LoRADS, cuLoRADS, can solve huge-scale semidefinite programming problems with remarkable efficiency. By effectively leveraging GPU computational power, cuLoRADS exhibits outstanding performance. Specifically, it can solve a set of MaxCut problems with $10^7 \times 10^7$ matrix variables in 10 seconds to 1 minute each on an NVIDIA H100 GPU with 80GB memory, whereas previous solvers demonstrated the capability of handling problems of this scale, required at least dozens of hours per problem on CPUs. Additionally, cuLoRADS shows exceptional scalability by solving 1) a MaxCut problem with a $170 \text{ million} \times 170 \text{ million}$ matrix variable and 2) a Matrix Completion problem with a 20 million $\times$ 20 million matrix variable and approximately 200 million constraints, both in a matter of minutes.

\end{abstract}
\noindent {\bf Keywords:}
large-/huge-scale SDP, GPU acceleration, low-rank factorization, logarithmic rank reduction

\section{Introduction}
Semidefinite programming (SDP) is a branch of convex optimization that generalizes linear programming to encompass semidefinite matrices. This paper addresses the following linear SDP problem:
\begin{equation}\label{prob:sdp}
\min_{X\in \mbb S^n} \inner{C}{X}\ \  \st \mcal A(X)=b,X\succeq 0,
\end{equation}
where $\mbb S^n$ denotes the space of $n \times n$ symmetric matrices, with $C \in \mbb S^n$. The linear map $\mcal A(\cdot):\mbb S^n\to \mbb R^m$ is specified as $\mcal A(X)=\bsbra{\inner{A_1}{X},\ldots,\inner{A_m}{X}}^\top$, where each $A_i$ belongs to $\mbb S^n$ for $i=1,\ldots,m$.

Semidefinite programming (SDP) has found extensive practical applications, including MaxCut~\citep{goemans1995improved}, optimal power flow~\citep{lavaei2011zero}, combinatorial optimization~\citep{boyd1997semidefinite}, sensor network localization~\citep{so2007theory}, and many others. Traditionally, interior point methods have been the go-to approach for solving SDPs due to their ability to achieve high accuracy in polynomial time~\citep{alizadeh1991combinatorial,alizadeh1995interior, nesterov1994interior,ye1997interior}. However, despite their effectiveness, interior point methods face scalability issues for large-scale problems, primarily due to the computational burden of repeatedly solving Newton systems. Recent research has shifted towards first-order methods to overcome these limitations and improve scalability ~\citep{wen2010alternating, lan2011primal, odonoghue2016conic}. \cite{zhao2010newton} proposed the SDPNAL algorithm, which employs the augmented Lagrangian method on the dual SDP and uses an inexact semi-smooth Newton-CG method to solve the inner subproblem. Subsequently, \cite{yang2015sdpnal+} improved SDPNAL by incorporating a majorized semi-smooth Newton-CG augmented Lagrangian method when solving degenerate SDPs. 

However, as problem sizes continue to increase, these methods inevitably suffer from the quadratically growing cost of maintaining the matrix variable. When the dimension reaches $n = 100,000$, storing a matrix variable of $n^2$ entries with double precision requires 75GB of memory, which may exceed the computational capacity of a single computer. To avoid processing such large/huge variables, several works~\citep{han2024low,monteiro2024low,yalccin2023semidefinite, BMImplement,marevcek2017low,erdogdu2022convergence} have explored finding low-rank solutions based on the Burer-Monteiro (BM) low-rank factorization of $X$, expressed as $X = RR^\top$, where matrix $R \in \mathbb{R}^{n \times r}$. 
In particular, \cite{monteiro2024low} recently proposed a new first-order method for solving large-scale SDPs with a bounded domain. The enhanced solver can solve the maximum stable set SDP with dimensions $(n,m) \approx (10^6, 10^7)$ within $10^{-5}$ relative precision. 
\cite{pataki1998rank} and \cite{barvinok1995problems} have shown that the problem admits an optimal low-rank solution with $r$ satisfying $r(r+1)/2 \leq m$. In addition, \cite{soandye2008} has theoretically demonstrated that employing a rank of $O(\log m)$ can efficiently produce an approximate solution of satisfactory quality. These results imply that choosing a lower rank $r$ can significantly reduce memory consumption without too much loss of accuracy. 

To enhance the performance of SDP solvers, we explore developing low-rank-based SDP algorithms on GPUs, which offer more efficient linear algebra operations than CPU architectures. 
Accelerating optimization solvers using GPUs has shown success in linear programming (LP)~\citep{lu2023cupdlp,lu2023cupdlpc} and quadratic programming (QP)~\citep{huang2024restarted,lu2023practical}. 
These solvers are based on primal-dual first-order methods, which can take advantage of GPUs' ability to perform matrix-vector multiplications efficiently. 
Therefore, a straightforward approach is to port existing SDP solvers, such as SDPLR~\citep{BMImplement}, HALLaR~\citep{monteiro2024low}, and LoRADS~\citep{han2024low}, to GPUs by utilizing GPU computation functions for matrix and vector operations. 
However, the inherent complexity of SDP computations, such as evaluating the linear mapping, presents challenges.  Directly rewriting a CPU solver for GPUs often leads to inefficiencies due to the different architecture designs.  
Our trials showed that simply porting the CPU implementation of LoRADS to GPU resulted in a solver that was even significantly slower than the original CPU version. The specific reasons for this will be discussed in Section~\ref{sec:low_rank_gpu}.

Thus, a natural question arises: \emph{Can we utilize the computational power of GPUs to achieve both efficiency and scalability in solving semidefinite programming problems?} In this work, we answer this question affirmatively and propose a solution to this challenge. Our main contributions are summarized as follows.

\begin{figure}[h]
    \centering
    \includegraphics[width=1\textwidth]{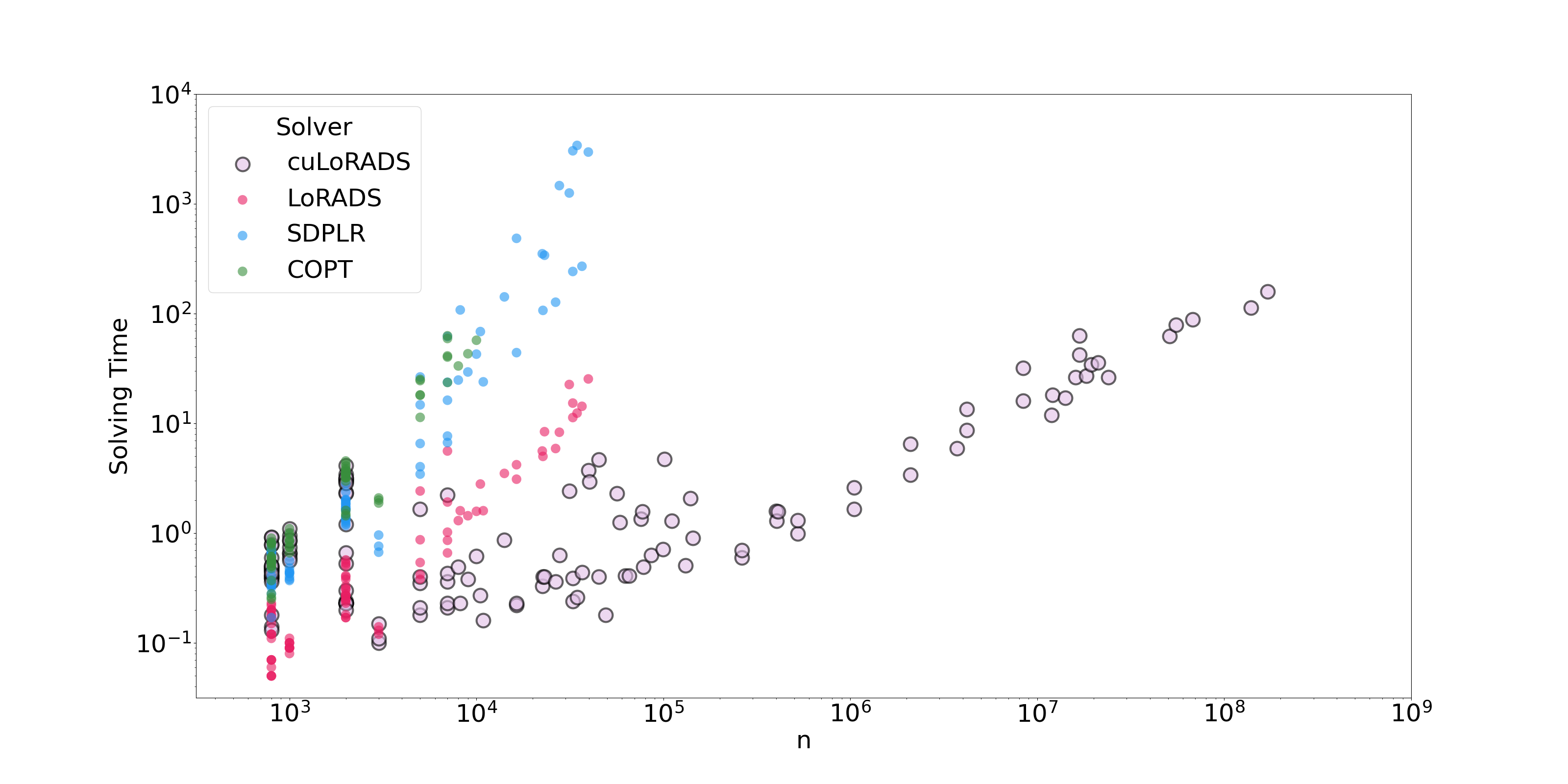}
    \caption{\footnotesize \textbf{Scalability and Efficiency of cuLoRADS on MaxCut Problems.} The figure shows the growth in solving times of the solvers as the problem dimensions increase (i.e., the dimension of the matrix variable $n$ and number of constraints $m$, where $m = n$ for MaxCut problems). The solving times and problem dimensions are presented on a logarithmic scale. The benchmark solvers involved are described in Section~\ref{sec:experimental_setting}. The figure shows that the solving times of cuLoRADS increase linearly with the problem dimensions.}
    \label{fig:maxcut}
\end{figure}

\begin{figure}[h]
    \centering
    \includegraphics[width=1\textwidth]{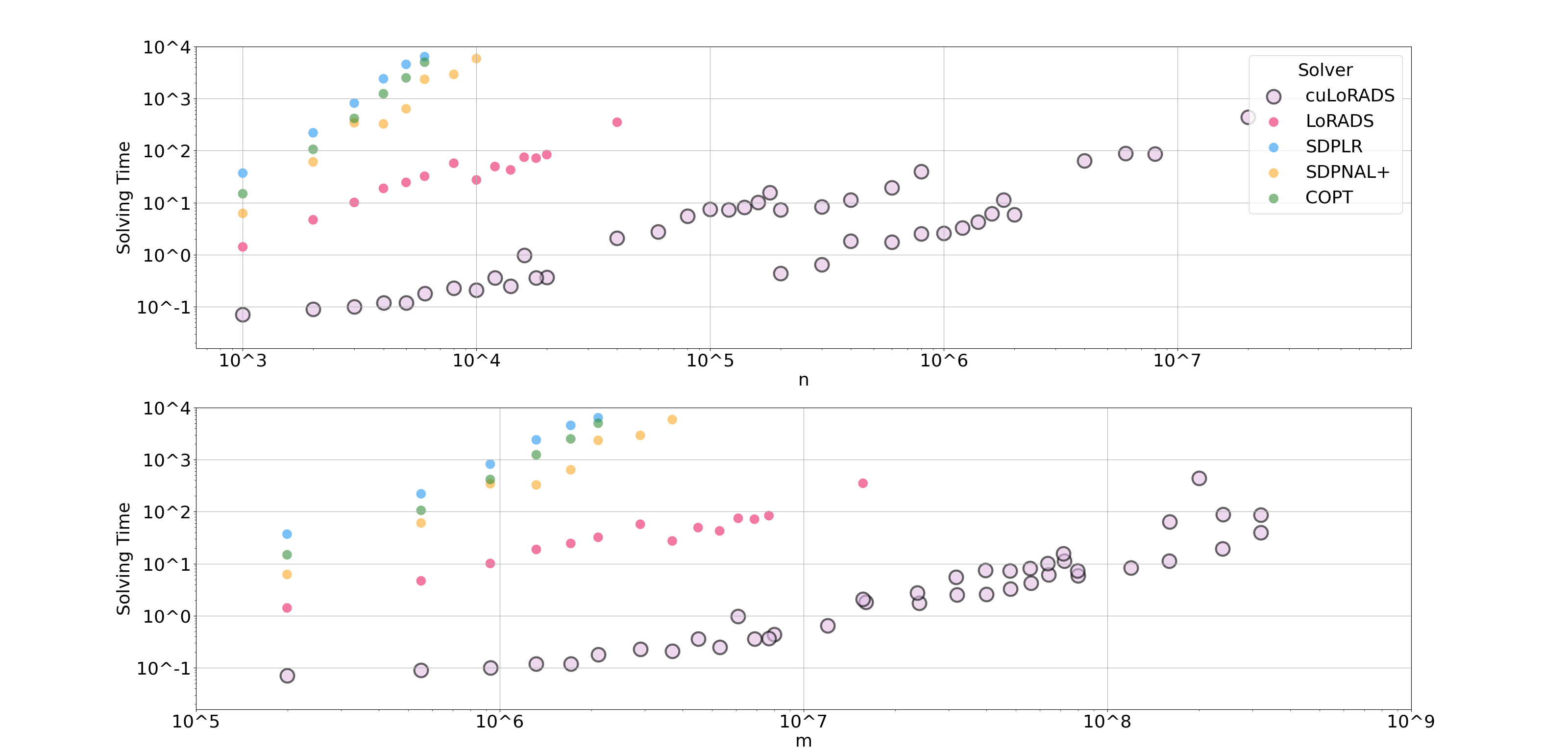}
    \caption{\footnotesize \textbf{Scalability and Efficiency of cuLoRADS on Matrix Completion Problems.} The figure shows the growth in solving times of the solvers as both the dimension of the matrix variable $n$ and the number of constraints $m$ increase. The solving times and problem dimensions are presented on a logarithmic scale. The figure shows that the solving times of cuLoRADS exhibit approximately linear growth as $n$ or $m$ increase.}
    \label{fig:mc}
\end{figure}

\begin{enumerate}
    \item 
    We propose accelerating a wide range of low-rank factorization methods for solving semidefinite programming problems by leveraging GPUs. Specifically, we effectively harness the computational power of GPUs to solve SDPs. We introduce a new, efficient method for calculating the vital operators in low-rank factorization methods on GPUs and propose a compression technique that exploits the column sparsity of the constraint matrix, reducing quadratic memory cost to the dependence on the number of nonzero elements of the problem and the size of low-rank matrix variables. In addition, we implement several new computation and storage patterns to reduce computation costs further and accelerate the overall process.
    We show that classical low-rank factorization-based methods, such as the one proposed by \cite{Burer2003ANP, BMImplement}, can be substantially improved in terms of scalability and efficiency for solving large-scale SDPs.
    \item We provide a GPU-based SDP solver, cuLoRADS, a GPU implementation of the Low-Rank ADMM Splitting Approach. cuLoRADS leverages both algorithmic and hardware advancements to achieve extraordinary solving speed and scalability.  Our empirical tests indicate that cuLoRADS significantly outperforms existing SDP solvers regarding running time on large-scale problems. Figure \ref{fig:maxcut} and Figure \ref{fig:mc} illustrate the efficiency and scalability of cuLoRADS on MaxCut and matrix completion problems. Remarkably, cuLoRADS is capable of solving a max-cut problem with a 170 million $\times$ 170 million matrix variable and 170 million constraints in under 160s. 
\end{enumerate}

\paragraph{Outline.} The paper is organized in the following manner. In the remainder of this section, we review the literature related to our work and the corresponding notations. In section~\ref{sec:cpu-sdp}, we introduce the CPU-based low-rank factorization SDP solvers, including the BM approach in Section~\ref{sec:bm} and LoRADS approach in Section~\ref{sec:lorads}. Section~\ref{sec:low_rank_gpu} presents a scalable and efficient approach for computing low-rank factorization methods on GPUs, using cuLoRADS as an example. In Section~\ref{sec:num}, we conduct comprehensive numerical experiments to evaluate the improvements brought about by our approach to classic low-rank factorization methods and to assess the performance of cuLoRADS. Finally, we conclude in Section~\ref{sec:conclusion}.


\subsection{Related Works}
Our paper is related to two streams of literature, including GPUs for mathematical optimization and low-rank SDP methods.
\paragraph{GPUs for mathematical optimization.}
The GPUs excel at efficiently performing matrix-vector multiplications, which means that first-order algorithms can seamlessly transition from CPUs to GPUs. However, first-order algorithms typically converge slowly, especially when transitioning from low- to high-precision, often requiring numerous iterations. The GPUs' ability to quickly execute matrix-vector multiplications significantly reduces the time per iteration, making GPU-based first-order algorithm solvers competitive with traditional solvers. The two GPU-based LP solvers (cuPDLP.jl~\citep{lu2023cupdlp} and cuPDLP-C~\citep{lu2023cupdlpc}) have been developed based on the PDLP algorithm~\citep{applegate2021practical}, which is a first-order method but have demonstrated impressive solver speed for LP problems. Moreover, two GPU-based first-order solvers for QP are proposed in~\cite{lu2023practical} and~\cite{huang2024restarted}. Given the success of the PDHG algorithm in LP, both papers adopted the PDHG framework to solve QP problems.  rAPDHG~\citep{lu2023practical} updates the primal variable directly by linearizing the objective function, while rPDHCG~\cite{huang2024restarted} employs the CG algorithm to achieve a more precise solution to the subproblem. Both algorithms are implemented on GPUs and demonstrate efficient solving performance on certain problems. It is worth noting that there are currently no GPU-based solvers specifically designed for large-/huge-scale SDPs.
\paragraph{The low-rank SDP methods.}
The standard SDPs often have low-rank solutions~\citep{pataki1998rank} or can be well approximated by low-rank matrices~\citep{soandye2008}. \cite{yurtsever2021scalable} aims to update the matrix variable using rank-one updates instead of updating a full-rank matrix each time. This approach of limiting the rank of the matrices to update similarly exploits the low-rank nature of the solutions. By focusing updates on certain low-rank matrix bases, extensive experiments have shown that this method can significantly accelerate the solving process. However, this method has a key drawback: it requires calculating eigenvalues to determine which low-rank update matrix is optimal. 
On the other hand, some works try to optimize the low-rank matrix directly~\citep{monteiro2024low,BMImplement,huang2024suboptimality}. They observe the low-rank solution $X=RR^\top$, so they chose to optimize $R$ directly rather than $X$. This decomposition method is also called the Burer-Monterio Low-Rank Factorization. Based on the above decomposition method, \cite{monteiro2024low,BMImplement} and~\cite{huang2024suboptimality} consider optimizing the following minimax optimization~\eqref{eq:bm_alm_func} using the ALM method.

\begin{equation}\label{eq:bm_alm_func}
    \min_{R\in\mbb R^{n\times r}}\max_{\lambda\in \mbb R^m}\mcal L_{\rho}(R,\lambda) = \inner{C}{RR^{\top}}+\inner{\mcal A(RR^{\top})-b}{\lambda} + \frac{\rho}{2}\norm{\mcal A(RR^{\top})-b}^2
\end{equation}
On the other hand, \cite{chen2023burer} and~\cite{han2024low} introduce an additional variable $(X=UV^\top, U=V)$ to increase the degrees of freedom of the variables. \cite{chen2023burer} consider a simple case where $r=1$. In this scenario, each iteration has an explicit solution.  In contrast, \cite{han2024low} considered a more general case and used a penalty method to incorporate the constraint $U=V$ into the objective function. They then solved minimax optimization~\eqref{eq:bilinear_admm_func} using the ADMM method, and the experimental results show good performance.

\begin{equation}\label{eq:bilinear_admm_func}
    \min_{U,V\in \mbb R^{n\times r}}\max_{\lambda\in \mbb R^m}\mcal L_{\rho}(U,V,\lambda) = \inner{C}{UV^{\top}}+\inner{\mcal A(UV^{\top})-b}{\lambda} + \frac{\rho}{2}\norm{\mcal A(UV^{\top})-b}^2 + \frac{\rho}{2}\norm{U-V}^2.
\end{equation}

\subsection{Notations}
Throughout the paper, we use $\norm{\cdot}_2$ to denote the Euclidean norm for vectors and the operator norm for matrices, and $\norm{\cdot}_F$ denotes the Frobenius norm for matrices. The inner product $\inner{\cdot}{\cdot}$ for two vectors is $\inner{x}{y}=x^\top y$ and that for two matrices is $\inner{U}{V}=trace(UV^\top)$. We use $\mcal A(X):=[\inner{A_1}{X},\ldots,\inner{A_m}{X}]^\top$ to denote the linear map. We denote the vectorization operation by $\vec(\cdot)$, where $\vec(X)=[x_1;x_2;\ldots;x_n]$ for $X=[x_1,\ldots,x_n]$. Then we can denote $\Abf = [\vec(A_1)^\top;\vec(A_2)^\top;\ldots;\vec(A_m)^\top]\in\mbb R^{m\times n^2}$ and rewrite the operator $\mcal A(X)$ as $\Abf\cdot \vec(X)$. For the adjoint operator $\mcal A^*(\lambda) = \sum_{i=1}^m\lambda_i A_i$, we can also rewrite it by $\mcal A^*(\lambda)=\mat(\Abf^\top \lambda)$, where $\mat(\cdot)$ is the inverse operator of $\vec(\cdot)$. Denote the minimum eigenvalue of the matrix $A$ as $\sigma_{\min}(A)$.

\section{CPU-based Low-Rank Factorization SDP Solvers\label{sec:cpu-sdp}}
 Recently, utilizing the property of low-rank solutions in SDP problems~\citep{soandye2008} can significantly reduce the dimension of optimization variables in extremely large-scale problems, making it possible to solve~\citep{BMImplement,han2024low}. SDPLR~\citep{BMImplement} use ALM to solve the SDP problem~\eqref{eq:bm_alm_func}, and the inner problem, the minimization of the quadratic problem, is solved by the limited-memory Broyden-Fletcher Goldfarb Shanno (LBFGS) method. The algorithm will not consume too much memory if this limited memory selection is sufficiently small. One of the most time-consuming operations in LBFGS is to calculate the gradient of~\eqref{eq:bm_alm_func} concerning $R$, i.e., $\nabla_{R}\mcal L_{\rho}(R,\lambda) = 2CR + 2\mcal A^*(\lambda)R + 2\rho\mcal A^*\brbra{\mcal A(RR^\top) - b}R,$ where $\mcal A^*(\lambda) = \sum_{i=1}^{m}\lambda_i A_i$ is the adjoint operator. 
 On the other hand, LoRADS~\citep{han2024low} uses a two-stage method to increase further the speed of solving SDPs. In the first stage, LoRADS uses SDPLR as a warm start and adopts ADMM to solve the SDP problem~\eqref{eq:bilinear_admm_func} in the second stage. The experiments show that ADMM converges faster than ALM(SDPLR) when variables enter some local regions. Every inner problem in LoRADS involving solving a linear system is addressed using the conjugate gradient method (CG)~\citep{nazareth2009conjugate}. Similarly to LBFGS, one of the most time-consuming operations in CG is to calculate the gradient of~\eqref{eq:bilinear_admm_func} concerning $U$ or $V$, i.e., $\nabla_U \mcal L_{\rho}(U, V,\lambda)= CV + \mcal A^*(\lambda)V + \rho \mcal A^*\brbra{\mcal A(UV^\top) - b}V + \rho(U-V)$ or $\nabla_{V}\mcal L_{\rho}(U, V,\lambda)$. In the section, we give a basic introduction to LoRADS in Section~\ref{sec:lorads}, which contains two parts: ALM warm start and ADMM acceleration. After that, we give a basic illustration of implementing LoRADS on CPUs, which helps us to understand the direction for improving our GPU solver. 

\subsection{The Burer-Monterio Approach\label{sec:bm}}
The low-rank decomposition $X=RR^\top$ incurs a natural algorithm that uses ALM to solve the low-rank SDP. The update of base algorithm ALM for solving~\eqref{eq:bm_alm_func} is 
\begin{equation}
\begin{cases}
R^{k+1} & =\argmin_{R}\mcal L_{\rho}(R,\lambda^{k})\\
\lambda^{k+1} & =\lambda^{k}+\rho\Brbra{\mcal A\brbra{R^{k+1}(R^{k+1})^{\top}} - b},
\end{cases}
\end{equation}
where $\rho$ is the penalty coefficient and the stepsize for updating dual variable $\lambda$. The update of dual variables is easy to understand. 
Here, we summarize the specific algorithm for updating $ R $ in SDPLR (the first stage of LoRADS) in Algorithm~\ref{alg:lbfgs}. The step size $\tau$ is chosen through an exact line search, detailed in Appendix \ref{sec:exact_line_search}.

\begin{algorithm}[ht]
\caption{LBFGS for updating $R$\label{alg:lbfgs}}
\begin{algorithmic}[1]
    \REQUIRE $T > 0$
    \FOR{$k=0,1,\ldots, K-1$}
        \STATE Set $D \leftarrow - \nabla_{R}\mathcal{L}_{\rho}(R,\lambda^k)$
        \FOR{$t=k-1,k-2,\cdots,k-T$}
            \STATE $\alpha_t=\beta_t \inner{s^t}{D}$, where $\beta_t= 1 / \inner{y^t}{s^t}$,
            \STATE $D \leftarrow D - \alpha_t y^t$
        \ENDFOR
        \FOR{$t=k-T,k-T+1,\cdots,k-1$}
            \STATE $D \leftarrow D + (\alpha_t-\beta_t \inner{y^t}{D})s^t$
        \ENDFOR
        \STATE $y^{k} = -\nabla_{R}\mathcal{L}_{\rho}(R^k,\lambda^k)$
        \STATE $R^{k+1} \leftarrow R^{k}+\tau \cdot D$
        \STATE Update $s^k = \tau \cdot D$, $y^k = y^k + \nabla_{R}\mathcal{L}(R^{k+1},\lambda^k)$, $\beta_k ={1}/{\inner{y^k}{s^k}}$
    \ENDFOR
    \STATE \textbf{Output: $R$}
\end{algorithmic}
\end{algorithm}

\subsection{The LoRADS Approach\label{sec:lorads}}

LoRADS~\citep{han2024low} is a step further from Burer-Monterio Approach. This method consists mainly of two parts: the first is optimizing~\eqref{eq:bm_alm_func} using the ALM method~\citep{BMImplement} (see Section~\ref{sec:bm}), and the second part involves using ADMM to optimize~\eqref{eq:bilinear_admm_func}.

In practice, transitioning from ALM to ADMM for solving~\eqref{eq:bilinear_admm_func} can expedite the process once ALM reduces the primal infeasibility to a specific threshold. The three steps of the base ADMM algorithm for solving~\eqref{eq:bilinear_admm_func} can be summarized as follows:
\begin{equation}
\begin{cases}
U^{k+1} & =\argmin_{U}\mcal L_{\rho}(U,V^{k},\lambda^{k})\\
V^{k+1} & = \argmin_{V}\mcal L_{\rho}(U^{k+1},V,\lambda^{k})\\
\lambda^{k+1} & =\lambda^{k}+\rho\Brbra{\mcal A\brbra{U^{k+1}(V^{k+1})^{\top}} - b}.
\end{cases}
\end{equation}
The core of ADMM for~\eqref{eq:bilinear_admm_func} is the update of $U$ and $V$. Indeed, the two minimization subproblems ($\nabla_{U}\mcal L_{\rho}(U, V^k,\lambda^k) = 0$) can be reformulated as two positive definite linear systems, which can be written as follows:
\begin{align}\label{eq:linsys}
    & \rho \mcal A^*\brbra{\mcal A(UV^\top)}V + \rho U=-CV-\mcal A^*(\lambda)V+\rho\mcal A^*(b)V + \rho V\\
    \iff &\Brbra{\sum_{i=1}^{m}\vec(A_i V)\vec(A_i V)^\top + \rho I}\vec(U)= \rho \vec(V)-\vec\brbra{CV+\mcal A^*(\lambda)V+\rho \mcal A^*(b)V}.\label{eq:apparent_linsys}
\end{align}
In practice, the linear system $\sum_{i=1}^{m}\vec (A_i V)\vec (A_i V)^\top + \rho I$ needs to store $nr\times nr$ elements, which is memory inefficient. On the other hand, the linear system is positive definite since $\rho>0$. Combining the two reasons, LoRADS uses CG to directly solve~\eqref{eq:linsys} but not~\eqref{eq:apparent_linsys} to avoid constructing the big linear system. For simplicity of notation, we denote $\mathfrak{A}(U) = \rho \mcal A^*\brbra{\mcal A(UV^\top)}V+\rho U, \mathfrak{b}=-CV-\mcal A^*(\lambda)V+\rho \mcal A^*(b)V+\rho V$. Now, we give the specific algorithm to update $U$ in Algorithm~\ref{alg:cg}. A similar procedure for updating $V$.
\begin{algorithm}[ht]
\caption{CG for updating $U$\label{alg:cg}}
\begin{algorithmic}[1]
    \REQUIRE $U_0,\vep>0,K>0$
    \STATE \textbf{Initialize:} $r_0=\mathfrak{b} - \mathfrak{A}(U_0)$
    \IF{$\|r_0\|_F\leq \vep$}
    \STATE{break}
    \ENDIF
    \STATE{$p_0=r_0, q_0=r_0$}
\FOR{$k=0,1,\ldots, K-1$}
\STATE{$Q=\mathfrak{A}(p_k)$}
\STATE{$\alpha_k=\tfrac{\inner{q_k}{r_k}}{\inner{p_k}{Q}}$}
\STATE{$U_{k+1}=U_k + \alpha_k p_k$}
\STATE{$r_{k+1}=r_k - \alpha_k Q$}
\IF{$ \|r_{k+1}\|_F \leq \vep$}
\STATE{break}
\ENDIF
\STATE{$q_{k+1}=r_{k+1}$}
\STATE{$\beta=\tfrac{\inner{q_{k+1}}{r_{k+1}}}{\inner{q_k}{r_k}}$}
\STATE{$p_{k+1}=r_{k+1}+\beta p_k$}
\ENDFOR
\STATE{\textbf{Output: $U$} }
\end{algorithmic}
\end{algorithm}

\subsection{Computation Patterns on CPUs}
Whether it's LBFGS or CG described in Section~\ref{sec:lorads}, aside from some unavoidable matrix inner products and matrix addition, the most critical aspect is how to compute $\mcal A(RR^\top)$ or $\mcal A(UV^\top), \mcal A^*(\lambda)$ and $\mcal A^*(\lambda)V$. The efficient computation of these three operations closely relates to how the data matrix A is stored.
Current CPU solvers, like open source solver LoRADS~\citep{han2024low}, SDPLR~\citep{BMImplement}, HDSDP~\citep{gao2022hdsdp}, mostly do the calculation for each small matrix $A_i$ separately. 
These solvers classify matrices based on certain characteristics, such as sparsity and low rank, and employ different computational methods customized for different categories of matrices to achieve optimal computational efficiency. 
On the other hand, the calculation order plays an important role in the process. Consider operation $\mcal A(UV^\top)$, we either calculate $X=UV^{\top}$ first and then calculate inner product $\inner{A_i}{X}$ for $m$ times or calculate $A_i U$ first and then calculate $\inner{(A_i U)}{V}$. The two methods have pros and cons. The first method allows calculation $UV^\top$ once, but it needs to store $n^2$ dense matrix $X$. The second one is memory efficiency when $r\ll n$, but the number of floating-point calculations it requires is greater than that of the first method.
Moreover, due to the different storage formats for the small matrices $A_i$ (sparse, dense, low-rank), the adjoint operation requires involving matrices of different formats for weighted summation. The calculation methods that can run efficiently on CPUs may not be suitable for efficient computation on GPUs. Lastly, the adjoint operation result $\mcal A^*(\lambda)$ is dense or sparse, which means we need to choose the proper matrix-matrix product method for this operation.

\section{Accelerating Low-Rank Factorization-Based Methods on GPU\label{sec:low_rank_gpu}} 
In this section, we discuss how to design an efficient and scalable GPU implementation of general first-order low-rank factorization-based methods, including the Burer-Monterio approach and LoRADS. Section \ref{sec:GPU} briefly introduces the hardware architecture and execution model of GPUs, highlighting the differences between CPUs and GPUs, and outlines general principles for developing optimization software with GPUs. Section~\ref{sec:basic} presents basic operations that are efficient on GPUs and describes an important structural design to ensure good performance. Sections~\ref{sec:efficiency} and~\ref{sec:scalability} discuss how to achieve efficiency and scalability on GPUs for low-rank factorization-based methods. Finally, in Section~\ref{sec:pattern}, we detail key adjustments and improvements in the computation and storage patterns on GPUs.

\subsection{CPU and GPU}\label{sec:GPU}

\paragraph{GPU overview}
Graphics Processing Units (GPUs) were initially developed for rendering graphics but have since evolved to handle a variety of general-purpose parallel computing tasks. GPUs are designed to maximize parallelism and throughput, achieved through the Single Instruction Multiple Data (SIMD) execution model. This model allows a single instruction to operate on multiple memory locations simultaneously and is supported by carefully optimized memory access patterns.

In Nvidia GPUs' operations, threads are organized into groups called warps, each typically consisting of 32 threads. Within a warp, all threads execute the same instruction at once but on different data pieces. This setup is especially efficient for parallel tasks like matrix and vector operations, allowing the GPU to use its numerous cores to execute a large number of computations simultaneously.

One crucial factor affecting GPU efficiency is memory access. Techniques such as memory coalescing~\citep{HwuProgramming} involve having adjacent threads access consecutive memory addresses. This helps to reduce the number of memory transactions and prevents cache pollution, leading to improved efficiency and throughput of GPU operations.

Kernel execution on a GPU is a process that begins with the CPU. The interaction between the CPU and GPU is facilitated through the PCIe (Peripheral Component Interconnect Express) bus. The CPU initiates kernel execution by using memory-mapped I/O to send execution commands to the GPU. These commands are written directly into a designated memory region that the GPU continuously monitors. Once the commands are received, the GPU executes the kernel code across its multiple cores, leveraging its SIMD execution model to perform parallel computations. 

\paragraph{Hardware-level differences between CPU and GPU. } In a nutshell, CPUs and GPUs are designed with fundamentally different goals: CPUs are latency-oriented, while GPUs are throughput-oriented~\citep{HwuProgramming}.
\begin{itemize}
\item CPUs are optimized to minimize latency for single-threaded tasks.  They achieve this through high clock speeds, deep and complex pipelines, and out-of-order execution~\citep{Tomasulo}, allowing instructions to be executed as soon as their operands are available, rather than strictly following their order in the code. This is complemented by a large, multi-level cache hierarchy that ensures fast access to frequently used data, thereby reducing latency~\citep{Pattersonquantitative}.

Additionally, CPUs employ sophisticated branch prediction~\citep{seznec2007256}
and speculative execution techniques ~\citep{speculationSteffan,speculationTorrellas}
to further minimize delays caused by conditional operations.
\item GPUs are optimized for throughput, focusing on maximizing the total volume of work performed within a given time frame through massive parallelism. This is facilitated by thousands of simpler, more power-efficient cores capable of handling numerous operations simultaneously. GPUs feature simpler pipelines and in-order execution to accommodate a greater number of cores per die. Instead of relying on large caches, GPUs utilize a high-bandwidth memory system that includes global memory and specialized memory types like shared and texture memory, which are essential for managing the voluminous data required for parallel processing.
\end{itemize}
These distinct design philosophies make CPUs ideal for tasks that require rapid, sequential execution and complex control logic, whereas GPUs are better suited for efficiently managing large-scale, data-parallel tasks.

\paragraph{CUDA}
CUDA (Compute Unified Device Architecture) is a parallel computing platform and application programming interface (API) model created by NVIDIA~\citep{HwuProgramming}. It allows developers to use familiar programming languages like C and C++ to write programs that execute on NVIDIA GPUs. CUDA provides direct access to the GPU's virtual instruction set and parallel computational elements, enabling developers to perform complex mathematical and computational tasks much more efficiently than on a CPU.

\paragraph{Enhancing computational power with cuBLAS and cuSPARSE}
To harness the advanced computational capabilities of GPUs, specialized libraries such as cuBLAS and cuSPARSE have been developed, targeting specific types of data and computations. cuBLAS, NVIDIA's implementation of the Basic Linear Algebra Subprograms (BLAS) for CUDA-enabled GPUs, offers highly optimized versions of standard BLAS functions. cuSPARSE delivers a comprehensive suite of sparse matrix algorithms optimized to exploit GPUs' unique memory access patterns and SIMD architecture. Both cuBLAS and cuSPARSE exemplify how specialized software can translate the raw power of GPU hardware into significant performance gains for complex mathematical computations.

\subsection{Basic Operations and Structral Design}\label{sec:basic}

\paragraph{Basic matrix and vector operations}
GPU computing is highly effective for simple and parallelizable tasks such as large matrix multiplications. It is well-recognized that matrix-matrix (MM), matrix-vector (MV), and vector-vector (VV) operations can be efficiently executed on GPUs using libraries such as cuBLAS and cuSPARSE. Therefore, low-rank factorization methods, including the Burer-Monterio approach and LoRADS, which rely on these operations, can potentially achieve significant acceleration with GPUs when properly designed.

\paragraph{The structural design}
To achieve reasonable performance on GPUs, it is crucial to carefully design the structural aspects of optimization software to minimize unnecessary data transfers. Taking cuLoRADS as an example, cuLoRADS executes major computations and iterations on the GPU, with only some pre- and post-processing steps performed on the CPU. This approach is designed to minimize data transfer during the solving process, which is vital given the expensive nature of CPU-GPU data communication in current GPU architectures. Although NVIDIA's new Blackwell architecture enhances CPU-GPU bandwidth and reduces data communication costs, economical management of data transfers remains essential for most existing GPU hardware.

\begin{figure}
    \centering
    \includegraphics[width=0.65\textwidth]{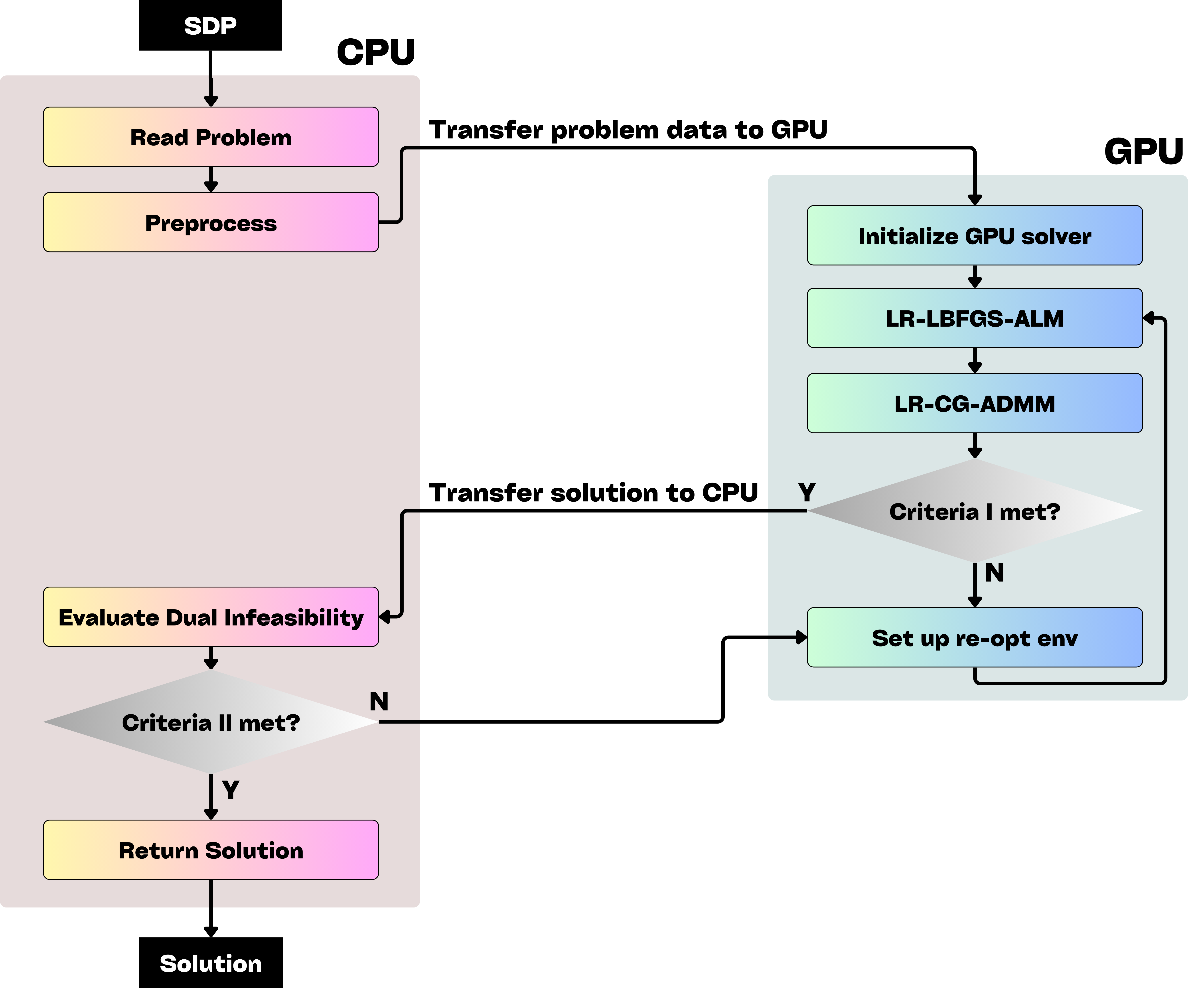}
    \caption{The internal workflow of cuLoRADS}
    \label{fig:flow_chart}
\end{figure}

As illustrated in Figure~\ref{fig:flow_chart}, data communication between the CPU and GPU occurs only during two scenarios: 1) transferring the loaded and preprocessed problem to the GPU and 2) transferring the solution back to the CPU. All algorithm iterations are conducted on the GPU, except for evaluating the dual infeasibility of the obtained solution.
This step involves calculating the minimum eigenvalue of a matrix, where we should employ first-order methods since we are focusing on large-scale problems.
Currently, there is no high-performance implementation of such first-order methods on GPU, so we temporarily leave this step to the CPU.

\subsection{Towards efficiency.}\label{sec:efficiency}
\begin{figure}
    \centering
    \includegraphics[width=0.8\textwidth]{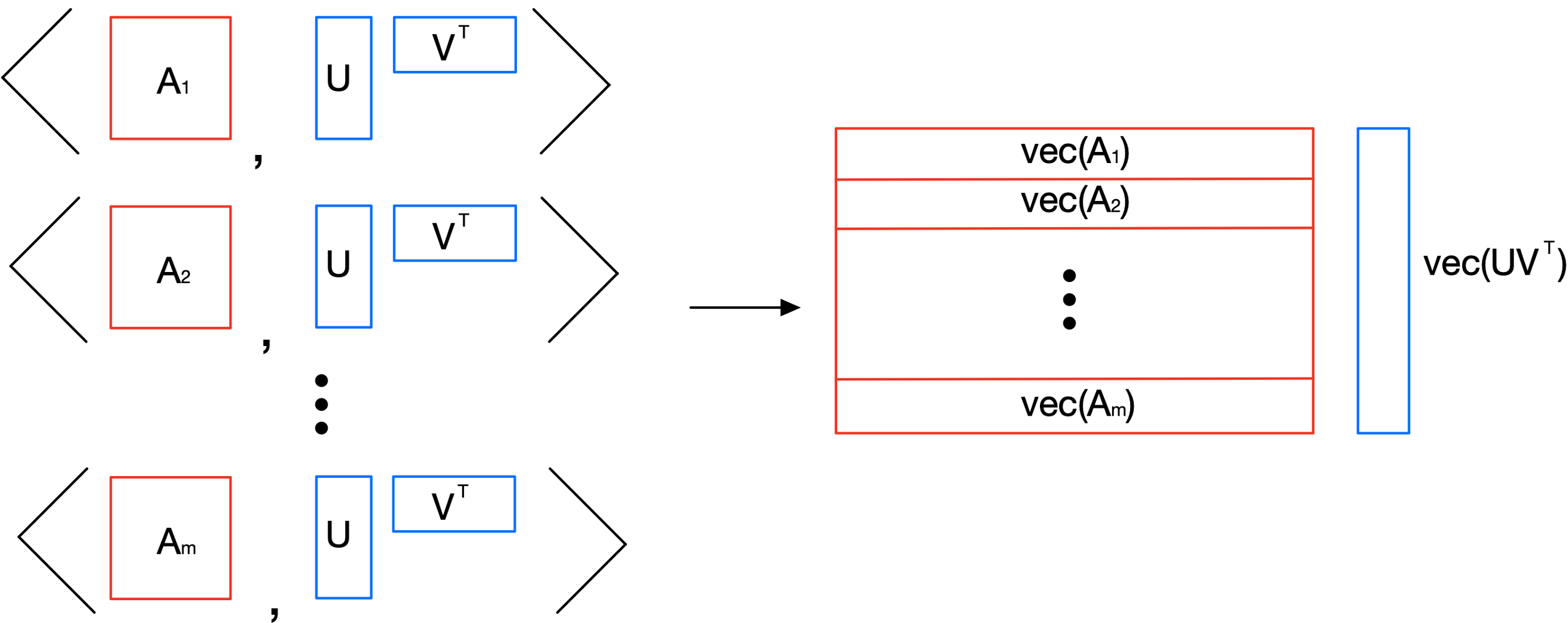}
    \caption{The modification of the linear mapping}
    \label{fig:AUV}
\end{figure}

This section discusses the bottleneck of efficiently calculating low-rank factorization-based methods on GPUs. While GPUs excel in large-scale, data-parallel tasks, numerous small kernels can sometimes lead to severe performance issues. 
The primary reason for this is that small kernels, without careful high-level parallelism design, fail to fully utilize a GPU's computing power, especially with advanced hardware featuring high computational capabilities. Self-implemented higher-level parallelism is typically less efficient than lower-level parallelism managed by the GPU driver.
Moreover, as discussed earlier, kernel launch overhead can become significant when there are numerous kernel calls. For example, consider the computation of $\mcal{A}(UV^\top)$. This involves calculating $X = UV^\top$ once and then computing the inner product $\langle A_i, X \rangle$ for $m$ times. Unfortunately, this requires $(m+1)$ matrix multiplication (MM) operations, resulting in $O(m)$ kernel launches. It's important to note that these kernel launches are inherently sequential, dictated by the scheduling and execution protocols of GPU kernels. As $m$ increases,  the cumulative overhead from these non-parallelizable kernel launches can become significant, potentially impacting the overall computational efficiency and scalability of the process. To illustrate, consider that each inner iteration of LR-LBFGS-ALM and each CG iteration of LR-CG-ADMM evaluates $\mcal{A}(RR^\top)$ or $\mcal{A}(UV^\top)$ $c_1$ and $c_2$ times, respectively. A common number for these iterations in LoRADS could be 1000. Even with a relatively small $m = 10000$, there would be $(1000 c_1 + 1000 c_2) \times (m+1) \geq 20$ million kernel launches. Assuming each kernel launch takes $10^{-5}$ seconds, the total overhead could exceed 200 seconds. This overhead renders the computation pattern impractical for GPUs, considering that solving a MaxCut problem with $n = m = 10000$ on a CPU using LoRADS takes only a few seconds.

Thus, the previously used computation pattern for the linear mapping on CPUs is unsuitable for GPUs. We introduce operation fusion for linear mapping to leverage the parallelism power of GPUs and avoid excessive kernel launches. Specifically, we horizontally vectorize each $A_i$ and concatenate them vertically into a large matrix $\mathbf{A}$, and vertically vectorize the result of $RR^\top$ or $UV^\top$. This allows us to combine $m$ inner product operations into a single matrix-vector multiplication $\mathbf{A} \text{vec}(RR^\top)$ or $\mathbf{A} \text{vec}(UV^\top)$. Figure \ref{fig:AUV} illustrates this concept. Following this approach, the computation of $\mcal{A}^*(\lambda)$ becomes $\mathbf{A}^\top \lambda$. With these modifications, we can efficiently compute the linear mappings on GPUs.

\subsection{Towards Scalability.}\label{sec:scalability}
\begin{figure}
    \centering
    \includegraphics[width=0.5\textwidth]{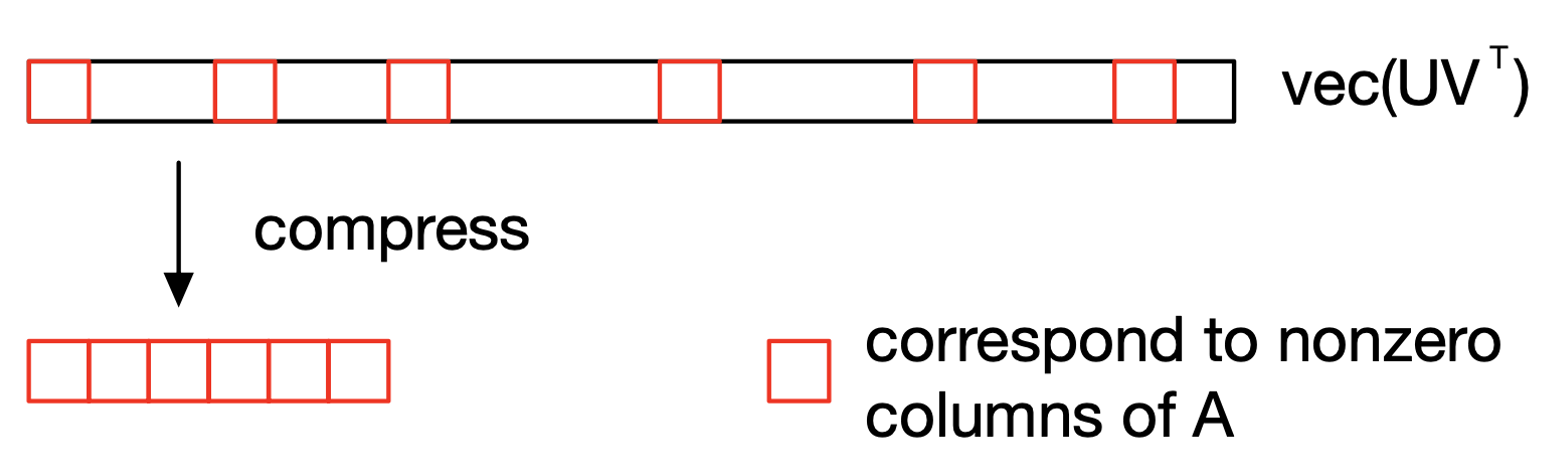}
    \caption{The modification of the linear mapping}
    \label{fig:compressing}
\end{figure}
In this section, we discuss the bottlenecks to scalability in solving large-scale SDP problems using low-rank factorization-based methods on GPUs. The operation fusion design above breaks the bottleneck to low-rank factorization-based methods' efficiency on GPUs. 
However, the scalability remains a problem. 
Note that $\text{vec}(RR^\top)$ is a dense vector with shape $n^2 \times 1$, and the storage needed for it is $O(n^2)$. 
And we can not avoid the $O(n^2)$ storage by using the "first $A_iU$ then $\inprod{(A_iU)}{V}$" pattern 
since we can not combine the $m$ operation into one under this computation sequence, which inevitably destroys the efficiency of the computation on GPU. 

\begin{figure}
    \centering
    \includegraphics[width=1\textwidth]{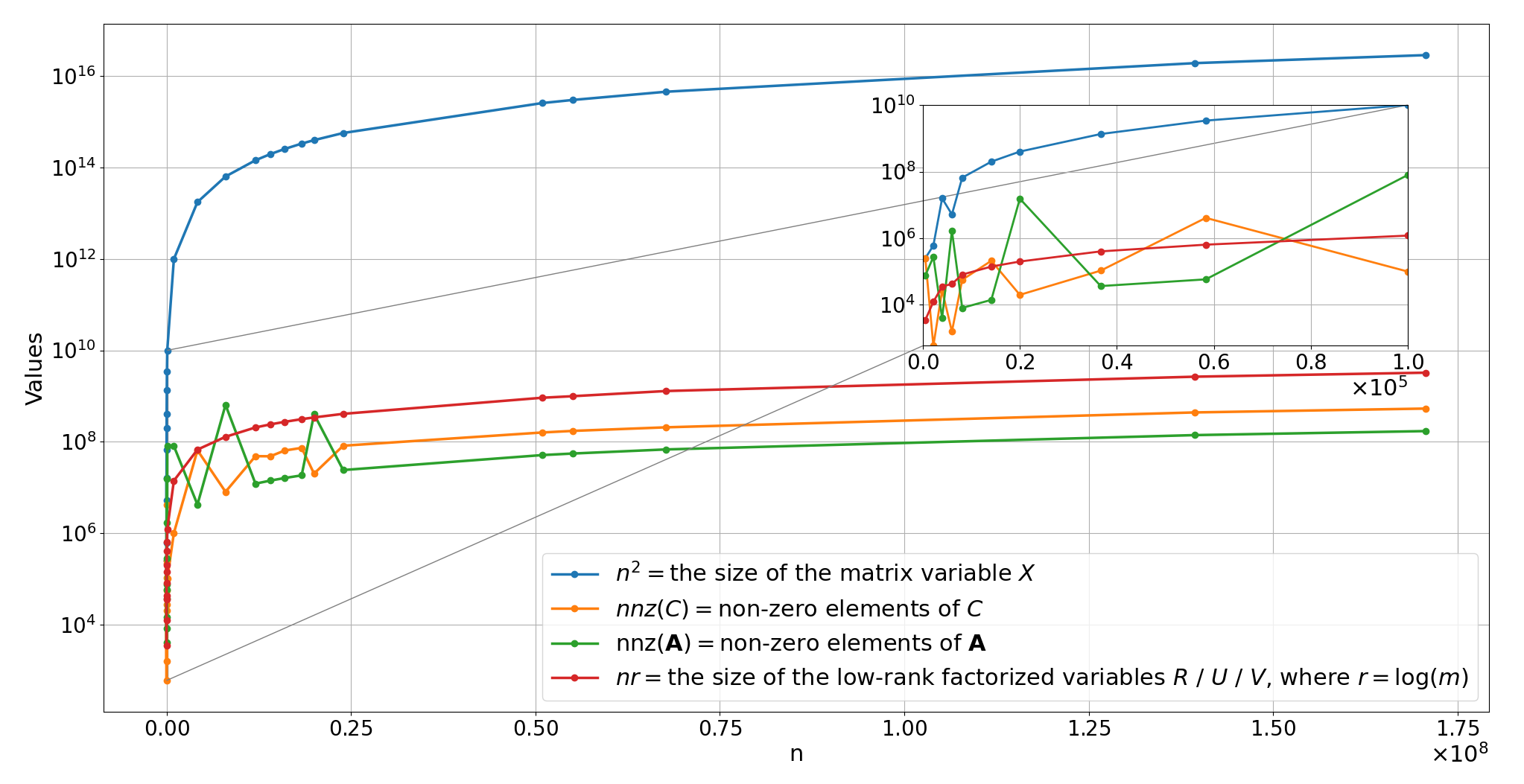}
    \caption{Comparison of Two Kinds of Memory Complexity}
    \label{fig:size}
\end{figure}

Suppose we are solving a SDP problem with the matrix variable $X \in \mathbb{R}^{n\times n}, n = 100000$, then the storage needed for $\text{vec}(RR^\top)$ or $\text{vec}(UV^\top)$ is about 75 GB, and 1863 GB when $X \in \mathbb{R}^{n\times n}, n=500000$. 
The VRAM of one of the currently most advanced GPUs, NVIDIA H100, is 80GB or 94GB. 
The quadratic growth of storage space needed makes it nearly impossible to solve huge-scale SDP problems on any modern GPU.
However, this does not mean that memory efficiency is impossible for LoRADS or other methods generally based on the low-rank factorization approach. 
Taking a look at the big constraint matrix $\mathbf{A}$, 
a potential fact is that the number of nonzero columns is less than or equal to the number of nonzero elements in $\mathbf{A}$.
And the only meaningful elements in $\text{vec}(RR^\top)$ or $\text{vec}(UV^\top)$ are those that correspond to the nonzero columns of $\mathbf{A}$ in position. 
Therefore, we can only retain the meaningful elements and ensure the storage needed for $\text{vec}(RR^\top)$ or $\text{vec}(UV^\top)$ would not exceed that of $\mathbf{A}$. 
This reduces the bottleneck for memory efficiency of the operation  from $O(n^2)$ to $O(\max(\nnz(\mathbf{A}), \nnz(C), nr))$, where $\nnz(\cdot)$ denotes the number of non-zeros elements in a matrix. 
This upper bound holds because the number of constraints is also less than or equal to $\nnz(\mathbf{A})$, so the storage needed for the dense vector $\lambda$ also would not exceed that of $\mathbf{A}$.
In real problems, the number of non-zeros columns of $\mathbf{A}$ can be far less than $\nnz(\mathbf{A})$, and $r$ can be approximately regarded as a constant if using logarithmic rank.
Figure \ref{fig:size} illustrates the difference between these two magnitudes as the dimension grows. When the problem dimension is large, the value of $\nnz(\mathbf{A})$, $\nnz(C)$, and $nr$ can be significantly smaller than $n^2$.

In practice, it needs some non-trivial adjustment to avoid the appearance of any $O(n^2)$ sized dense data. Denote the number of columns of $\mathbf{A}$ by $\text{nnzcol}\mathbf(A)$, 
and the compressed vector of $\text{vec}(UV^\top)$ by $\bar{x} \in \mathbb{R}^{\text{nnzcol}\mathbf(A) \times 1}$. We need to obtain $\bar{x}$ directly from $U$ and $V$ without creating any $O(n^2)$ sized dense data. 
No known package by CUDA or third party can conduct such an operation.
We implemented a customized kernel called CompressedOuterProduct() for this operation, which takes in $U$ and $V$ (or $R$) and returns $\bar{x}$. Within CompressedOuterProduct(), each thread is responsible for calculating one element in $\bar{x}$, i.e. $x_k = U_{i, \cdot} V_{j, \cdot}^\top$.
To avoid index mapping in each thread of every kernel calling, we pre-calculate two index mapping vectors for directly locating the corresponding row in $U$ and column in $V^\top$ needed for calculating the specific element in $\bar{x}$.

\subsection{Some Important Computation and Storage Patterns} \label{sec:pattern}

Finally, we discuss some computational and storage patterns crucial for GPU solvers. For large-scale SDP problems, storing the objective matrix $C$ in a sparse format is essential. Otherwise, it would require storing a dense $n \times n$ matrix. By applying the compression technique introduced above, we align the shape of $\mathbf{A}$ with $\bar{x}$ to simplify the implementation of the matrix-vector multiplication $\mathbf{A} \bar{x}$ and enhance the operation's efficiency. During preprocessing, we delete the zero columns of $\mathbf{A}$. Although this deletion does not reduce the storage used by $\mathbf{A}$ in sparse form, it accelerates calculations involving CUSPARSE on $\mathbf{A}$ by leveraging the shorter column index range after column compression.

The storage design for $\mathbf{A}^\top$ is slightly more complex. Considering the operation $CV + \mcal A^*(\lambda)V + \rho \mcal A^*(b)V$ in \eqref{eq:linsys}, it can be written as:
\begin{equation} \label{eq:At}
    \text{mat}\brbra{\text{vec}(C) + \mathbf{A}^\top \lambda + \rho \mathbf{A}^\top b}V, 
\end{equation}
where $\text{mat}(\cdot)$ denotes reshaping an $n^2 \times 1$ vector into an $n \times n$ square matrix. If we delete the zero rows of $\mathbf{A}^\top$ as we do for $\mathbf{A}$, the dimensions of $\mathbf{A}^\top \lambda$ and $\text{vec}(C)$ would not match, nor would the dimensions of $\mathbf{A}^\top b$ and $\text{vec}(C)$. To align these dimensions, we could wrap $\mathbf{A}^\top \lambda$ and $\mathbf{A}^\top b$ with $n^2 \times 1$ sparse vectors. However, this would lead to inefficient sparse-add-sparse operations, requiring the reconstruction of nonzero indices.

To optimize the operation in \eqref{eq:At}, we maintain a set of indices $\Omega := \Omega_{C} \cup \Omega_{\mathbf{A}}$, where $\Omega_{C} := \{i \mid i \in [n^2], \text{vec}(C)_i \neq 0\}$ and $\Omega_{\mathbf{A}} := \{j \mid j \in [n^2], \mathbf{A}_{i,j} = 0, \forall i \in [m]\}$. The indices in $\Omega$ correspond to the union of nonzero indices in the vectorized objective matrix $C$ and the positions of nonzero columns in $\mathbf{A}$. We compress $\mathbf{A}^\top$ so that the remaining rows match the elements in $\Omega$. At the beginning of the solving process, we construct a compressed version of $\text{vec}(C)$ that stores the elements of $C$ in $\Omega$ as a dense vector $\mathfrak{c}$. Consequently, any dot product of $\mathbf{A}^\top$ and a vector has the same shape and one-to-one index correspondence with $\mathfrak{c}$. The additions in \eqref{eq:At} become $\mathfrak{c} + \mathbf{A}^\top \lambda + \rho \mathbf{A}^\top b$. Finally, we wrap the result in an $n \times n$ sparse matrix to conduct the dot product with $V$.

Additionally, using the sparse format for $C$ is less efficient when $C$ is dense and small. We store $C$ as a dense matrix in such scenarios and use dense operators for all related computations. In this situation, neither the columns of $\mathbf{A}$ nor the rows of $\mathbf{A}^\top$ are compressed. The usage of sparse or dense storage for $C$ is automatically decided in cuLoRADS.

\section{Numerical Study\label{sec:num}}
In this section, we demonstrate the effectiveness of our approach by examining the numerical performance of cuLoRADS. We compare cuLoRADS with its CPU counterpart, LoRADS, implemented in C programming language \citep{han2024low}. We also present results from other SDP solvers, including SDPLR, SDPNAL+, and COPT, on several not-that-large-scale problems. Additionally, we test cuLoRADS on a range of large- to huge-scale problems that can only be solved by cuLoRADS. Section \ref{sec:experimental_setting} details the experimental settings. In Section \ref{sec:maxcut_results}, we evaluate cuLoRADS against the aforementioned solvers on MaxCut problems. Section~\ref{sec:matrix_completion_results} demonstrates the results of the solvers in solving Matrix Completion problems. Finally, Section~\ref{sec:benchmark_results} tests cuLoRADS against the other four solvers on a broader range of problems from Hans Mittelmann's SDP benchmarks and SDPLIB.
\subsection{Implementation Details}
\paragraph{LoRADS.} The LoRADS solver is implemented with the C programming language. It incorporates several important techniques and heuristic strategies that enhance its performance, including the dynamic rank strategy, the control of the penalty factor $\rho$, and the heuristic factor. For details on these techniques, refer to Sections 4.3 and 5.1 in~\cite{han2024low}.

\paragraph{cuLoRADS.} The cuLoRADS solver is implemented with the Julia programming language\footnote{The cuLoRADS solver is available at \url{https://github.com/COPT-Public/cuLoRADS}.}. It inherits the techniques applied in LoRADS. Additionally, to address the issue where LoRADS sometimes fails to balance the feasibility and optimality of the solution—specifically, the trade-off between primal infeasibility, dual infeasibility, and the primal-dual gap—cuLoRADS incorporates a re-optimization (\texttt{re-opt}) technique. This technique improves the dual infeasibility and primal-dual gap errors when necessary through rescaling, thereby ensuring balanced convergence of all three criteria: primal feasibility, dual feasibility, and the primal-dual gap.

The idea behind the \reopt technique is that rescaling the objective matrix $C$ can balance the feasibility and optimality conditions. Through extensive experimentation, we observed the following phenomena. When scaling down the objective matrix, i.e., multiplying it by a factor less than 1, the emphasis on optimality increases, resulting in a higher-quality solution at the same level of primal feasibility. Conversely, scaling up the objective matrix, i.e., multiplying it by a factor greater than 1, can relax the optimality requirements, allowing the solver to meet the primal infeasibility tolerance faster, albeit with a solution of lower optimality quality.

Before further discussion, it is important to note that repeatedly evaluating dual infeasibility is impractical for large-scale problems, as this step involves calculating the smallest eigenvalue of an $n \times n$ matrix. Therefore, we only evaluate primal infeasibility and the primal-dual gap during the solving process and assess final dual infeasibility only when we are confident in the current primal-dual solution pair, i.e., when both the primal infeasibility and primal-dual gap tolerances are met.

The \reopt technique can be detailed as follows. Initially, we stop the optimization process once the primal infeasibility tolerance is met, as discussed in Section 5.2 of~\cite{han2024low}. However, for some problems, due to the imbalance of errors, dual infeasibility and the primal-dual gap may not be satisfied. We enter the \reopt phase after satisfying the primal infeasibility tolerance while other errors remain unmet. Essentially, the \reopt phase involves solving a scaled-down problem again using the previous output solution as a warm start. Upon stopping, we re-evaluate the errors; if the final criteria are not met, we repeat the \reopt process until the final criteria are satisfied or the maximum number of \reopt iterations is reached.

In our experiments, the \reopt factor, i.e., the factor by which the objective matrix is scaled, is 0.1. We provide two levels of \reopt: one that stops when the mild final criteria (primal feasibility and primal-dual gap) are met and another that stops when all three criteria (primal infeasibility, dual infeasibility, and primal-dual gap) are met. The \reopt level is controlled by the parameter \texttt{reoptLevel} in cuLoRADS, with \texttt{reoptLevel} = 1 corresponding to the mild \reopt level and \texttt{reoptLevel} = 2 requiring all three errors to be met. For large-scale problems, the mild \reopt level (\texttt{reoptLevel} = 1) is often preferred, as evaluating dual infeasibility can be very time-consuming.

\paragraph{The Burer-Monteiro method.} The method proposed by \cite{Burer2003ANP} is incorporated into both LoRADS and cuLoRADS. Specifically, we follow the implementation of its enhanced version by \cite{BMImplement}, applying all the aforementioned techniques in LoRADS and cuLoRADS.

\subsection{Experimental Setting\label{sec:experimental_setting}}
In this section, we introduce the experimental setting, including the solvers considered, the errors evaluated, and the stopping criteria.

\paragraph{Solvers.} We test the following solver in our numerical experiment.
\begin{itemize}
\item cuLoRADS. The Julia implementation of LoRADS on GPU.
\item LoRADS~\citep{han2024low}. The C implementation of LoRADS on CPU.
\item SDPLR~\citep{BMImplement}. An open-source solver that employs the classic Burer-Monteiro low-rank factorization method based on the ALM, which uses the limited-memory BFGS algorithm and an exact line search procedure.
\item SDPNAL+~\citep{sun2019sdpnal}. A MATLAB solver for large-scale semidefinite programming problems, utilizing a majorized semi-smooth Newton-CG augmented Lagrangian method. 
\item COPT~\citep{ge2022cardinal}. A commercial solver utilizing the interior point method, COPT is currently acknowledged as the top-performing solver in Hans Mittelmann's benchmarks for SDP problems\footnote{\url{https://plato.asu.edu/ftp/sparse_sdp.html}}.
\end{itemize}

\paragraph{Stopping criteria and errors reported.}
Each solver has its embedded stopping criteria, some of which cannot be modified. The default stopping criteria for the solvers are as follows:

\begin{itemize}
\item cuLoRADS terminates differently depending on the $\texttt{reoptLevel}$ setting:

If $\texttt{reoptLevel} = 0$ (no \texttt{re-opt}), cuLoRADS terminates when 
$$
\frac{\left\|\mathcal{A}\left(X\right) - b \right\|_2}{1+\left\|b\right\|_\infty} \leq \epsilon.
$$

If $\texttt{reoptLevel} = 1$ (mild \texttt{re-opt}), cuLoRADS terminates when 
$$
\max \left\{\frac{\norm{\mcal A(X)-b}_2}{1+\norm{b}_1}, \frac{\inner{C}{X}-\lambda^{\top}b}{1+\abs{\inner{C}{X}}+\abs{\lambda^{\top}b}} \right\} < \epsilon.
$$

If $\texttt{reoptLevel} = 2$ (strict \texttt{re-opt}), cuLoRADS terminates when 
$$
\max \left\{\frac{\norm{\mcal A(X)-b}_2}{1+\norm{b}_1},\frac{\abs{\min\{0,\sigma_{\min}(C-\mcal A^{*}(\lambda))\}}}{1+\norm{c}_1}, \frac{\inner{C}{X}-\lambda^{\top}b}{1+\abs{\inner{C}{X}}+\abs{\lambda^{\top}b}} \right\} < \epsilon.
$$

\item LoRADS and SDPLR terminate when 
$$
\frac{\left\|\mathcal{A}\left(X\right) - b \right\|_2}{1+\left\|b\right\|_\infty} \leq \epsilon.
$$

\item SDPNAL+ terminates when 

\begin{equation*}
 \max \left\{
\frac{\|\mathcal{A}(X)-b\|}{1+\|b\|},\frac{\left\|\mathcal{A}^{*}(y)+S-C\right\|}{1+\|C\|},\frac{1}{5} \frac{\left\|X-\Pi_{S_+^{n}}(X-S)\right\|}{1+\|X\|+\|S\|}\right\}\leq \epsilon,
\end{equation*}
where $S$ is the slack variable of the dual of \eqref{prob:sdp}, and $\Pi_{S_+^{n}}(X-S)$ is the metric projection of $(X-S)$ onto $S_+^{n}$.

\item COPT terminates when 
$$
\max \left\{\frac{\norm{\mcal A(X)-b}_2}{1+\norm{b}_1},\frac{\abs{\min\{0,\sigma_{\min}(C-\mcal A^{*}(\lambda))\}}}{1+\norm{\vec(C)}_1}, \frac{\inner{C}{X}-\lambda^{\top}b}{1+\abs{\inner{C}{X}}+\abs{\lambda^{\top}b}} \right\} < \epsilon.
$$ 
\end{itemize}

We set $\epsilon = 10^{-5}$ for all the experiments in the following section. Note that sometimes primal infeasibility is normalized with the infinity norm $\|b\|_\infty$ rather than $\|b\|_1$. As the number of constraints grows, primal infeasibility normalized with $\|b\|_\infty$ can be much more rigorous than that normalized with $\|b\|_1$.

To evaluate the quality of the solutions produced by the solvers in a unified manner, we consider the following three commonly used errors:

Primal Infeasibility:
$$\text{err}_1 = \frac{\norm{\mcal A(X)-b}_2}{1+\norm{b}_1},$$

Dual Infeasibility:
$$\text{err}_2 = \frac{\abs{\min\{0,\sigma_{\min}(C-\mcal A^{*}(\lambda))\}}}{1+\norm{\vec(C)}_1},$$

Primal-Dual Gap:
$$\text{err}_3 = \frac{\inner{C}{X}-\lambda^{\top}b}{1+\abs{\inner{C}{X}}+\abs{\lambda^{\top}b}}.$$

Additionally, a time limit of 10,000 seconds is set for each experiment.

\paragraph{Hardware.} For GPU benchmarks, we use an NVIDIA H100 with 80 GB VRAM deployed on a cluster with an Intel Xeon Platinum 8469C CPU. The CPU usage of Julia is restricted to 1 thread. For CPU benchmarks, we use a MacBook Pro with a 12-core M3 Pro chip and 18 GB RAM\footnote{In the following experiments, for the CPU solvers tested on the same problems as those in \cite{han2024low}, we directly use the results from that work since the settings and hardware are identical.}.

\subsection{The MaxCut Problems\label{sec:maxcut_results}}
The MaxCut problem in a weighted graph $ G $ with edge weights $ w $ involves finding a division of its vertices into two distinct groups such that the sum of the weights of the edges between these groups is maximized. A semidefinite programming approach can be utilized to solve this problem, formulated as follows:
\begin{align*}
    \max_{X \in {\mathbb{S}^{n \times n}}} \ \ &\langle \frac{1}{4}L(G,w), X \rangle \\
    \text{s.t.} \ \ &X_{ii} = 1 \ \ , \forall i \in [n] \\
    &X \succeq 0,
\end{align*}
where $ L(G,w) $ is the Laplacian matrix of the graph, defined by:
\begin{align*}
L(G, w)_{ij}:= \begin{cases} 
-w_{ij} & \text{if } (i, j) \in E \\ 
\sum_{k} w_{ik} & \text{if } i = j \\ 
0 & \text{otherwise}
\end{cases}.
\end{align*}

We consider MaxCut problems from two sources, one is the Gset\footnote{See \url{https://web.stanford.edu/~yyye/yyye/Gset/}.}~\citep{ye_gset_2003}; and also a set of large- to huge-scale MaxCut problems generated from the sparse matrices in the University of Florida Sparse Matrix Collection\footnote{See \url{https://sparse.tamu.edu/}.}~\citep{sparse_mat_set}. SDPNAL+ is not included for MaxCut problems as it is inefficient for this type of problem, according to \cite{wang2023solving, han2024low}.

\paragraph{The Gset.} The Gset includes MaxCut problems ranging in size from $n = m = 800$ to $n = m = 10,000$. Tables \ref{tab:gset1} and \ref{tab:gset2} show the performance of all five solvers. Table~\ref{tab:gset1} includes the smaller problems with $n = m \le 3,000$, where cuLoRADS, LoRADS, SDPLR, and COPT have similar performance. However, on the larger problems in Table \ref{tab:gset2}, cuLoRADS and LoRADS start to show their advantages, with cuLoRADS being slightly faster. From the two tables, we can see that cuLoRADS exhibits competitive performance compared to its CPU counterpart, LoRADS, on these small to moderate-sized MaxCut problems.

\begin{table}[!ht]
    \centering
    \scriptsize
    \caption{Numerical Results on Gset}
    \vspace{5pt}
    \setlength{\tabcolsep}{1.2mm}
    \begin{tabular}{|cc|cc|cc|cc|cc|cc|}
    \hline
        ~ & ~ & \multicolumn{2}{c|}{\textbf{cuLoRADS}} & \multicolumn{2}{c|}{\textbf{LoRADS}} & \multicolumn{2}{c|}{\textbf{SDPLR}} & \multicolumn{2}{c|}{\textbf{COPT}} \\ 
        name & n & time & $\text{error}_{\text{max}}$ & time & $\text{error}_{\text{max}}$ & time & $\text{error}_{\text{max}}$ & time &  $\text{error}_{\text{max}}$ \\ \hline
        G1 & 800 & 0.4 & 2.4E-07 & 0.2 & 2.5E-07 & 0.8 & 1.6E-06 & 0.8 & 2.4E-09 \\ 
        G2 & 800 & 0.4 & 2.2E-07 & 0.2 & 2.4E-07 & 0.7 & 1.8E-06 & 0.8 & 3.2E-09 \\ 
        G3 & 800 & 0.4 & 3.8E-07 & 0.2 & 5.6E-07 & 0.8 & 3.8E-07 & 0.9 & 3.9E-09 \\ 
        G4 & 800 & 0.4 & 1.7E-07 & 0.2 & 4.6E-07 & 0.7 & 4.1E-07 & 0.8 & 2.3E-09 \\ 
        G5 & 800 & 0.4 & 2.0E-07 & 0.2 & 3.2E-07 & 0.7 & 2.6E-07 & 0.7 & 3.4E-09 \\ 
        G6 & 800 & 0.6 & 4.3E-07 & 0.2 & 4.6E-07 & 0.7 & 9.8E-08 & 0.6 & 2.4E-08 \\ 
        G7 & 800 & 0.5 & 5.6E-08 & 0.2 & 1.6E-07 & 0.6 & 6.4E-06 & 0.6 & 1.2E-08 \\ 
        G8 & 800 & 0.5 & 4.7E-07 & 0.2 & 6.6E-07 & 0.7 & 3.6E-06 & 0.6 & 1.1E-08 \\ 
        G9 & 800 & 0.5 & 3.7E-08 & 0.2 & 1.0E-07 & 0.5 & 1.5E-06 & 0.5 & 1.2E-08 \\ 
        G10 & 800 & 0.5 & 5.7E-08 & 0.2 & 8.4E-08 & 0.6 & 1.9E-06 & 0.8 & 9.4E-09 \\ 
        G11 & 800 & 0.2 & 1.0E-07 & 0.1 & 3.0E-07 & 0.5 & 9.3E-07 & 0.3 & 6.4E-08 \\ 
        G12 & 800 & 0.1 & 2.4E-07 & 0.1 & 4.2E-07 & 0.3 & 9.6E-07 & 0.3 & 2.8E-08 \\ 
        G13 & 800 & 0.1 & 2.0E-07 & 0.1 & 1.6E-07 & 0.2 & 3.3E-06 & 0.3 & 4.1E-08 \\ 
        G14 & 800 & 0.5 & 9.6E-07 & 0.1 & 2.0E-06 & 0.3 & 3.7E-07 & 0.5 & 2.5E-08 \\ 
        G15 & 800 & 0.5 & 7.3E-07 & 0.1 & 1.2E-06 & 0.3 & 1.4E-06 & 0.6 & 6.1E-09 \\ 
        G16 & 800 & 0.5 & 1.1E-06 & 0.1 & 1.0E-06 & 0.3 & 1.2E-06 & 0.6 & 1.0E-08 \\ 
        G17 & 800 & 0.4 & 4.2E-07 & 0.1 & 1.2E-06 & 0.3 & 3.5E-06 & 0.5 & 7.5E-09 \\ 
        G18 & 800 & 0.9 & 1.9E-07 & 0.1 & 3.0E-07 & 0.4 & 2.0E-06 & 0.4 & 3.4E-08 \\ 
        G19 & 800 & 0.8 & 4.2E-07 & 0.1 & 2.6E-07 & 0.4 & 4.6E-06 & 0.5 & 6.3E-08 \\ 
        G20 & 800 & 0.8 & 2.4E-07 & 0.1 & 3.7E-07 & 0.3 & 1.5E-06 & 0.5 & 4.7E-08 \\ 
        G21 & 800 & 0.9 & 1.9E-07 & 0.1 & 1.7E-07 & 0.4 & 5.9E-07 & 0.5 & 4.1E-08 \\ 
        G43 & 1000 & 0.7 & 3.3E-07 & 0.1 & 3.7E-07 & 0.4 & 1.1E-06 & 1.0 & 4.8E-09 \\ 
        G44 & 1000 & 0.7 & 4.8E-07 & 0.1 & 4.8E-07 & 0.4 & 6.3E-07 & 0.7 & 7.5E-09 \\ 
        G45 & 1000 & 0.6 & 2.5E-07 & 0.1 & 2.3E-07 & 0.4 & 1.3E-06 & 1.0 & 1.8E-09 \\ 
        G46 & 1000 & 0.7 & 2.7E-07 & 0.1 & 3.0E-07 & 0.6 & 9.2E-07 & 0.9 & 5.9E-09 \\ 
        G47 & 1000 & 0.6 & 2.7E-07 & 0.1 & 2.9E-07 & 0.4 & 1.3E-06 & 1.1 & 8.1E-09 \\ 
        G51 & 1000 & 0.9 & 8.8E-07 & 0.1 & 9.0E-07 & 0.5 & 9.3E-07 & 0.8 & 7.5E-09 \\ 
        G52 & 1000 & 1.0 & 7.5E-07 & 0.1 & 7.4E-07 & 0.4 & 3.9E-07 & 0.9 & 4.6E-09 \\ 
        G53 & 1000 & 1.1 & 2.2E-07 & 0.1 & 6.0E-07 & 0.4 & 3.3E-06 & 0.8 & 4.8E-09 \\ 
        G54 & 1000 & 0.9 & 5.9E-07 & 0.1 & 6.0E-07 & 0.5 & 9.0E-07 & 0.7 & 9.4E-09 \\ 
        G22 & 2000 & 3.0 & 8.7E-07 & 0.3 & 1.0E-06 & 1.7 & 9.6E-07 & 3.4 & 7.1E-09 \\ 
        G23 & 2000 & 2.9 & 2.9E-08 & 0.3 & 2.9E-08 & 1.3 & 1.9E-06 & 3.7 & 4.2E-08 \\ 
        G24 & 2000 & 2.4 & 3.0E-08 & 0.3 & 2.8E-08 & 1.9 & 2.0E-07 & 3.5 & 4.5E-08 \\ 
        G25 & 2000 & 2.3 & 4.1E-08 & 0.3 & 1.6E-08 & 1.3 & 9.0E-07 & 3.3 & 1.1E-08 \\ 
        G26 & 2000 & 3.1 & 4.6E-07 & 0.3 & 6.2E-07 & 1.7 & 2.2E-06 & 3.8 & 6.5E-09 \\ 
        G27 & 2000 & 3.2 & 5.4E-07 & 0.4 & 5.4E-07 & 1.7 & 5.7E-07 & 3.5 & 7.6E-09 \\ 
        G28 & 2000 & 3.4 & 4.7E-07 & 0.5 & 7.5E-07 & 2.0 & 6.4E-06 & 4.1 & 5.6E-09 \\ 
        G29 & 2000 & 3.2 & 3.7E-07 & 0.3 & 4.4E-07 & 2.4 & 3.9E-07 & 4.3 & 4.9E-09 \\ 
        G30 & 2000 & 2.9 & 1.4E-07 & 0.3 & 5.8E-08 & 1.4 & 1.6E-06 & 4.5 & 6.8E-09 \\ 
        G31 & 2000 & 4.1 & 1.4E-06 & 0.4 & 1.1E-07 & 2.0 & 5.8E-07 & 3.9 & 6.8E-09 \\ 
        G32 & 2000 & 0.2 & 2.1E-07 & 0.2 & 1.9E-07 & 1.4 & 1.6E-06 & 1.5 & 2.9E-08 \\ 
        G33 & 2000 & 0.2 & 1.4E-07 & 0.2 & 5.2E-08 & 1.2 & 1.7E-06 & 1.4 & 7.0E-09 \\ 
        G34 & 2000 & 0.2 & 1.2E-07 & 0.2 & 2.9E-07 & 1.2 & 1.0E-06 & 1.6 & 8.7E-09 \\ 
        G35 & 2000 & 0.2 & 1.9E-07 & 0.2 & 2.9E-07 & 1.6 & 5.5E-07 & 3.5 & 8.2E-09 \\ 
        G36 & 2000 & 0.2 & 5.7E-07 & 0.2 & 3.1E-07 & 1.4 & 1.4E-06 & 3.0 & 6.1E-09 \\ 
        G37 & 2000 & 0.2 & 1.8E-07 & 0.3 & 2.4E-07 & 1.8 & 4.5E-07 & 3.5 & 7.0E-09 \\ 
        G38 & 2000 & 0.2 & 7.0E-07 & 0.3 & 2.7E-07 & 1.9 & 4.3E-07 & 3.3 & 3.9E-09 \\ 
        G39 & 2000 & 0.3 & 4.8E-09 & 0.4 & 1.5E-08 & 1.8 & 1.9E-07 & 3.3 & 5.3E-09 \\ 
        G40 & 2000 & 0.5 & 9.0E-09 & 0.6 & 1.3E-08 & 2.8 & 5.4E-07 & 3.5 & 5.3E-09 \\ 
        G41 & 2000 & 0.7 & 6.4E-09 & 0.5 & 5.6E-08 & 2.0 & 2.3E-07 & 3.2 & 1.3E-08 \\ 
        G42 & 2000 & 1.2 & 7.4E-07 & 0.5 & 1.0E-06 & 1.7 & 6.8E-07 & 3.2 & 9.0E-09 \\ 
        G48 & 3000 & 0.1 & 4.5E-07 & 0.1 & 1.2E-06 & 0.7 & 1.7E-06 & 2.1 & 9.6E-09 \\ 
        G49 & 3000 & 0.1 & 5.9E-07 & 0.1 & 5.1E-07 & 0.8 & 1.0E-06 & 2.0 & 4.7E-09 \\ 
        G50 & 3000 & 0.2 & 6.3E-07 & 0.1 & 1.1E-06 & 1.0 & 4.9E-06 & 1.9 & 5.0E-09 \\ 
 \hline
    \end{tabular}
    \label{tab:gset1}
\end{table}

\begin{table}[!ht]
    \centering
    \scriptsize
    \caption{Numerical Results on Gset, Cont}
    \vspace{5pt}
    \setlength{\tabcolsep}{1.2mm}
    \begin{tabular}{|cc|cc|cc|cc|cc|}
    \hline
        ~ & ~ & \multicolumn{2}{c|}{\textbf{cuLoRADS}} & \multicolumn{2}{c|}{\textbf{LoRADS}} & \multicolumn{2}{c|}{\textbf{SDPLR}} & \multicolumn{2}{c|}{\textbf{COPT}} \\ 
        name & n & time & $\text{error}_{\text{max}}$ & time & $\text{error}_{\text{max}}$ & time & $\text{error}_{\text{max}}$ & time &  $\text{error}_{\text{max}}$ \\ \hline
        G55 & 5000 & 0.2 & 9.2E-08 & 0.4 & 1.4E-07 & 3.5 & 4.3E-07 & 18.0 & 7.2E-09 \\ 
        G56 & 5000 & 0.2 & 2.0E-07 & 0.4 & 2.5E-07 & 4.0 & 8.5E-07 & 18.2 & 5.0E-09 \\ 
        G57 & 5000 & 0.4 & 1.8E-08 & 0.5 & 9.4E-09 & 6.6 & 3.2E-07 & 11.3 & 1.1E-08 \\ 
        G58 & 5000 & 0.4 & 2.1E-07 & 0.9 & 3.3E-07 & 14.8 & 3.7E-06 & 24.5 & 6.5E-09 \\ 
        G59 & 5000 & 1.7 & 5.0E-08 & 2.4 & 1.2E-07 & 26.5 & 1.0E-06 & 25.2 & 8.4E-09 \\ 
        G60 & 7000 & 0.2 & 1.7E-07 & 0.7 & 1.6E-07 & 7.7 & 6.4E-07 & 40.2 & 4.2E-09 \\ 
        G61 & 7000 & 0.2 & 2.2E-07 & 0.9 & 2.3E-07 & 6.7 & 2.2E-06 & 41.2 & 8.5E-09 \\ 
        G62 & 7000 & 0.4 & 4.9E-09 & 1.0 & 4.4E-08 & 16.3 & 1.4E-07 & 23.6 & 6.8E-09 \\ 
        G63 & 7000 & 0.4 & 2.1E-07 & 1.9 & 2.0E-07 & 23.7 & 3.7E-07 & 59.3 & 4.0E-09 \\ 
        G64 & 7000 & 2.2 & 2.9E-08 & 5.6 & 1.9E-07 & 61.8 & 6.0E-07 & 62.8 & 3.0E-09 \\ 
        G65 & 8000 & 0.5 & 5.4E-09 & 1.3 & 9.1E-09 & 24.8 & 5.4E-07 & 33.4 & 5.6E-09 \\ 
        G66 & 9000 & 0.4 & 6.9E-09 & 1.4 & 1.0E-08 & 29.4 & 5.3E-08 & 43.1 & 9.6E-09 \\ 
        G67 & 10000 & 0.6 & 5.9E-08 & 1.6 & 2.0E-08 & 42.7 & 9.3E-07 & 57.2 & 8.3E-09 \\ \hline
    \end{tabular}
    \label{tab:gset2}
\end{table}

\paragraph{Large- to huge-scale MaxCut problems.}
We select a set of random sparse matrices from the University of Florida Sparse Matrix Collection, specifically from the DIMACS10, SNAP, and GenBank groups. Since we are dealing with large-scale problems, we only include low-rank solvers in our analysis. Considering the hardware's storage limitations and the software's memory efficiency, we present the results of LoRADS and SDPLR for MaxCut problems with $n = m \le 40,000$ in Table \ref{tab:maxcut_l}. Table \ref{tab:maxcut_l} shows the results of cuLoRADS on larger MaxCut problems with dimensions up to $10^8$.

In Table \ref{tab:maxcut_l}, cuLoRADS begins to demonstrate its efficiency compared to CPU solvers. We also observe the good scalability of cuLoRADS, as the solving time does not notably increase when the problem dimension grows by a factor of 10 (e.g., solving cit-HepPh with $n = 27,770$ and delaunay\_n18 with $n = 262,111$ both in 0.6 seconds). Furthermore, cuLoRADS can solve MaxCut problems with $n \approx 10^7$ to an accuracy of $\max (\text{err}_1, \text{err}_2, \text{err}_3) \lesssim 10^{-8}$ in 10 seconds to 1 minute. According to the work by \cite{yurtsever2021scalable}, SketchyCGAL—the only other solver reported to handle SDP problems of this scale—requires about 30-50 hours to solve each problem to an accuracy where the three errors are under $10^{-1}$ (noting that they use $\|\cdot \|_2$ to normalize primal infeasibility and dual infeasibility). The largest problem, with $n \approx 1.7 \times 10^8$, is solved by cuLoRADS in 158 seconds.

\begin{table}[!ht]
    \centering
    \scriptsize
    \caption{Numerical Results on Large-Scale MaxCut Problems}
    \vspace{5pt}
    \setlength{\tabcolsep}{1.2mm}
    \begin{tabular}{|cc|cc|cc|cc|}
    \hline
        ~ & ~ & \multicolumn{2}{c|}{\textbf{cuLoRADS}} & \multicolumn{2}{c|}{\textbf{LoRADS}} & \multicolumn{2}{c|}{\textbf{SDPLR}} \\ 
        name & n & time & $\text{error}_{\text{max}}$ & time & $\text{error}_{\text{max}}$ & time & $\text{error}_{\text{max}}$ \\ \hline
        delaunay\_n13 & 8192 & 0.2 & 3.5E-06 & 1.6 & 3.6E-06 & 108 & 4.5E-08 \\ 
        vsp\_p0291\_seymourl\_iiasa & 10498 & 0.3 & 7.8E-06 & 2.8 & 6.7E-06 & 68.6 & 1.3E-07 \\ 
        p2p-Gnutella04 & 10879 & 0.2 & 7.8E-06 & 1.6 & 2.3E-05 & 23.9 & 7.1E-07 \\ 
        vsp\_befref\_fxm\_2\_4\_air02 & 14109 & 0.9 & 1.8E-06 & 3.5 & 4.0E-06 & 142 & 5.8E-08 \\ 
        delaunay\_n14 & 16384 & 0.2 & 3.0E-06 & 4.2 & 3.6E-06 & 485 & 2.0E-08 \\ 
        fe\_sphere & 16386 & 0.2 & 2.5E-06 & 3.1 & 9.8E-06 & 44.1 & 2.2E-07 \\ 
        cs4 & 22499 & 0.3 & 4.8E-06 & 5.6 & 1.4E-06 & 351 & 1.2E-07 \\ 
        p2p-Gnutella25 & 22687 & 0.4 & 4.1E-06 & 5.0 & 3.0E-06 & 107 & 4.8E-08 \\ 
        ca-CondMat & 23133 & 0.4 & 5.5E-06 & 8.4 & 6.9E-06 & 340 & 6.7E-08 \\ 
        p2p-Gnutella24 & 26518 & 0.4 & 3.1E-06 & 5.9 & 8.4E-06 & 127 & 7.2E-07 \\ 
        cit-HepTh & 27770 & 0.6 & 3.2E-06 & 8.3 & 3.2E-06 & 1466 & 2.5E-08 \\ 
        as-caida & 31379 & 2.4 & 1.0E-05 & 22.6 & 1.6E-05 & 1255 & 3.2E-07 \\ 
        delaunay\_n15 & 32768 & 0.4 & 9.5E-07 & 15.3 & 1.9E-06 & 3042 & 3.9E-07 \\ 
        se & 32768 & 0.2 & 2.6E-07 & 11.3 & 4.0E-07 & 242 & 4.7E-08 \\ 
        cit-HepPh & 34546 & 0.3 & 4.6E-06 & 12.4 & 4.7E-06 & 3407 & 1.1E-07 \\ 
        p2p-Gnutella30 & 36682 & 0.4 & 1.0E-06 & 14.3 & 1.4E-06 & 270 & 1.7E-07 \\ 
        vsp\_south31\_slptsk & 39668 & 3.7 & 3.9E-06 & 25.4 & 6.6E-06 & 2964 & 4.6E-07 \\ \hline
    \end{tabular}
    \label{maxcut_l}
\end{table}

\begin{table}[!ht]
    \centering
    \scriptsize
    \caption{Numerical Results on Large- to Huge-Scale MaxCut Problems}
    \vspace{5pt}
    \setlength{\tabcolsep}{1.2mm}
    \begin{tabular}{|cc|cc||cc|cc|}
    \hline
        ~ & ~ & \multicolumn{2}{c||}{\textbf{cuLoRADS}} & ~ & ~ & \multicolumn{2}{c|}{\textbf{cuLoRADS}} \\ 
        name & n & time & $\text{error}_{\text{max}}$ & name & n & time & $\text{error}_{\text{max}}$ \\ \hline
        vsp\_sctap1-2b\_and\_seymourl & 40174 & 3.0  & 3.5E-06 & delaunay\_n20 & 1048576 & 1.7  & 1.1E-07 \\ 
        fe\_body & 45087 & 0.4  & 1.3E-06 & rgg\_n\_2\_20\_s0 & 1048576 & 2.6  & 6.2E-08 \\ 
        vsp\_model1\_crew1\_cr42\_south31 & 45101 & 4.7  & 3.9E-06 & delaunay\_n21 & 2097152 & 3.4  & 9.8E-08 \\ 
        cca & 49152 & 0.2  & 1.1E-09 & rgg\_n\_2\_21\_s0 & 2097152 & 6.5  & 6.6E-08 \\ 
        vsp\_bump2\_e18\_aa01\_model1\_crew1 & 56438 & 1.9  & 6.3E-07 & 333SP & 3712815 & 5.9  & 4.8E-08 \\ 
        vsp\_bcsstk30\_500sep\_10in\_1Kout & 58348 & 1.3  & 3.9E-07 & delaunay\_n22 & 4194304 & 8.7  & 6.3E-08 \\ 
        p2p-Gnutella31 & 62586 & 0.4  & 1.1E-06 & rgg\_n\_2\_22\_s0 & 4194304 & 13.5  & 3.4E-08 \\ 
        delaunay\_n16 & 65536 & 0.4  & 1.4E-06 & delaunay\_n23 & 8388608 & 16.0  & 5.7E-08 \\ 
        soc-Epinions1 & 75888 & 1.4  & 2.7E-06 & rgg\_n\_2\_23\_s0 & 8388608 & 32.0  & 6.3E-09 \\ 
        vsp\_vibrobox\_scagr7-2c\_rlfddd & 77328 & 1.6  & 4.9E-07 & asia\_osm & 11950757 & 12.0  & 1.9E-08 \\ 
        fe\_tooth & 78136 & 0.5  & 5.2E-07 & hugetrace-00010 & 12057441 & 18.1  & 6.0E-09 \\ 
        vsp\_c-60\_data\_cti\_cs4 & 85830 & 0.6  & 8.4E-07 & road\_central & 14081816 & 17.0  & 7.0E-09 \\ 
        fe\_rotor & 99617 & 0.7  & 3.6E-07 & hugetrace-00020 & 16002413 & 26.2  & 1.5E-08 \\ 
        vsp\_mod2\_pgp2\_slptsk & 101364 & 0.9  & 4.9E-06 & delaunay\_n24 & 16777216 & 42.1  & 2.3E-08 \\ 
        598a & 110971 & 1.3  & 3.3E-07 & rgg\_n\_2\_24\_s0 & 16777216 & 62.7  & 1.7E-08 \\ 
        delaunay\_n17 & 131072 & 0.5  & 7.4E-07 & hugebubbles-00000 & 18318143 & 27.1  & 3.5E-09 \\ 
        vsp\_finan512\_scagr7-2c\_rlfddd & 139752 & 2.1  & 2.8E-07 & hugebubbles-00010 & 19458087 & 34.2  & 5.5E-09 \\ 
        fe\_ocean & 143437 & 0.9  & 2.2E-07 & hugebubbles-00020 & 21198119 & 35.9  & 3.9E-09 \\ 
        amazon0302 & 262111 & 0.6  & 2.3E-07 & road\_usa & 23947347 & 26.2  & 2.9E-09 \\ 
        delaunay\_n18 & 262144 & 0.7  & 3.1E-07 & europe\_osm & 50912018 & 62.1  & 2.2E-09 \\ 
        amazon0312 & 400727 & 1.6  & 2.2E-07 & kmer\_V2a & 55042369 & 79.1  & 7.5E-09 \\ 
        amazon0601 & 403394 & 1.3  & 1.8E-07 & kmer\_U1a & 67716231 & 88.4  & 4.9E-09 \\ 
        amazon0505 & 410236 & 1.6  & 1.8E-07 & kmer\_P1a & 139353211 & 114  & 8.7E-09 \\ 
        delaunay\_n19 & 524288 & 1.0  & 1.3E-07 & kmer\_A2a & 170728175 & 158  & 1.1E-08 \\ 
        rgg\_n\_2\_19\_s0 & 524288 & 1.3  & 6.6E-08 & ~ & ~ & ~ & ~ \\ \hline
    \end{tabular}
    \label{tab:maxcut_l}
\end{table}
\subsection{The Matrix Completion Problems\label{sec:matrix_completion_results}}
Consider a matrix $M \in \mathbb{R}^{n \times m}$ with observed entries in $\Omega \subseteq \{1, \ldots, n \} \times \{ 1, \ldots, m \}$. The goal of the matrix completion problem is to determine a matrix $Y$ such that $Y_{ij} = M_{ij}$ for all $(i,j) \in \Omega$, while minimizing the rank of $Y$. This problem can be expressed as the following convex optimization problem:

\begin{align*}
\min_{Y \in \mathbb{R}^{m \times n}} \ \ \norm{Y}_* \ \ 
\text{subject to} \ \  Y_{ij} = M_{ij}, \ \ \forall (i,j) \in \Omega,
\end{align*}
where $\norm{Y}_*$ represents the nuclear norm of $Y$. Additionally, the nuclear norm minimization problem can be reformulated into a semidefinite programming problem as follows:

\begin{align*}
\min_{Y \in \mathbb{R}^{m \times n}} \ \ \text{trace}(W_1) + \text{trace}(W_2) \ \ 
\text{subject to} \ \  Y_{ij} = M_{ij}, \ \ \forall (i,j) \in \Omega, 
\begin{bmatrix}
W_1 & Y^{\top} \\
Y & W_2
\end{bmatrix} \succeq 0,
\end{align*}
which is equivalent to:
\begin{align*}
\min_{X \in \mathbb{S}^{(m+n) \times (m+n)}} \ \ \langle I, X \rangle \ \ 
\text{subject to} \ \ \left\langle \begin{bmatrix}
0_{m \times m} & E_{ij}^{\top} \\
E_{ij} & 0_{n \times n}
\end{bmatrix}, X \right\rangle = 2 M_{ij}, \ \ \forall (i,j) \in \Omega, \ X \succeq 0,
\end{align*}
where $X = \begin{bmatrix}
W_1 & Y^{\top} \\
Y & W_2
\end{bmatrix}$, and $E_{ij}$ is the $ij$-th adjacency matrix with $(E_{ij})_{ij} = 1$ and zeros elsewhere.

Table \ref{tab:mc_l} compares the solvers' performance on matrix completion problems. CuLoRADS and LoRADS are significantly more efficient than other solvers, with cuLoRADS being much more efficient than LoRADS. cuLoRADS achieves more than 100 times speedup on a problem with a $40,000 \times 40,000$ matrix variable and approximately 16 million constraints.

Table \ref{tab:mc_h} demonstrates the performance of cuLoRADS on larger problems that none of the other solvers could handle on the used CPU hardware. Compared to its CPU counterpart, cuLoRADS solves a problem with an $8 \text{million} \times 8 \text{million}$ matrix variable and 320 million constraints in about 89 seconds, while LoRADS can only solve a problem with a $20,000 \times 20,000$ matrix variable and 8 million constraints in a similar time, improving the scale of the problem solved by hundreds of times within the same solving time. The largest matrix completion problem solved by cuLoRADS features a $20 \text{million} \times 20 \text{million}$ matrix variable and 200 million constraints, completed in less than 8 minutes.

\begin{table}[!ht]
    \centering
    \small
    \caption{Numerical Results on Matrix Completion Problems}
    \vspace{5pt}
    \begin{tabular}{|cc|cc|cc|cc|cc|cc|}
    \hline
        ~ & ~ & \multicolumn{2}{c|}{\textbf{cuLoRADS}} & \multicolumn{2}{c|}{\textbf{LoRADS}} & \multicolumn{2}{c|}{\textbf{SDPLR}} & \multicolumn{2}{c|}{\textbf{SDPNAL+}} & \multicolumn{2}{c|}{\textbf{COPT}} \\ 
        n & m & time & $\text{error}_{\text{max}}$ & time & $\text{error}_{\text{max}}$ & time & $\text{error}_{\text{max}}$ & time &  $\text{error}_{\text{max}}$ & time &  $\text{error}_{\text{max}}$ \\ \hline 
        1000 & 199424 & 0.1  & 4.9E-06 & 1.4  & 3.8E-06 & 37.0  & 2.26E-05 & 6.2  & 1.5E-06 & 14.8  & 1.3E-06 \\ 
        2000 & 550536 & 0.1  & 2.2E-06 & 4.7  & 1.9E-06 & 220 & 2.39E-06 & 60.7  & 2.0E-06 & 106 & 2.0E-06 \\ 
        3000 & 930328 & 0.1  & 1.7E-06 & 10.1  & 9.3E-07 & 819 & 1.60E-05 & 343 & 1.5E-06 & 416 & 4.8E-06 \\ 
        4000 & 1318563 & 0.1  & 5.3E-07 & 18.8  & 5.8E-06 & 2401 & 4.07E-06 & 326 & 3.0E-06 & 1239 & 6.8E-06 \\ 
        5000 & 1711980 & 0.1  & 3.9E-06 & 24.5  & 1.2E-06 & 4547 & 1.88E-06 & 636 & 2.1E-07 & 2502 & 6.7E-06 \\ 
        6000 & 2107303 & 0.2  & 3.0E-07 & 32.2  & 2.2E-06 & 6390 & 1.05E-05 & 2347 & 1.8E-06 & 5006 & 5.1E-06 \\ 
        8000 & 2900179 & 0.2  & 1.5E-06 & 57.2  & 7.3E-07 & t & - & 2921 & 6.4E-06 & t & - \\ 
        10000 & 3695929 & 0.2  & 7.0E-06 & 27.3  & 1.2E-06 & t & - & 5879 & 1.4E-05 & t & - \\ 
        12000 & 4493420 & 0.4  & 2.6E-07 & 49.4  & 2.3E-06 & t & - & t & - & t & - \\ 
        14000 & 5291481 & 0.3  & 1.1E-06 & 42.6  & 1.5E-06 & t & - & t & - & t & - \\ 
        16000 & 6089963 & 1.0  & 4.5E-06 & 74.7  & 5.2E-06 & t & - & t & - & t & - \\ 
        18000 & 6889768 & 0.4  & 2.8E-07 & 71.6  & 4.3E-06 & t & - & t & - & t & - \\ 
        20000 & 7688309 & 0.4  & 7.4E-06 & 83.8  & 9.5E-07 & t & - & t & - & t & - \\ 
        40000 & 15684167 & 2.1  & 1.2E-06 & 351 & 1.4E-06 & t & - & t & - & t & - \\ \hline
    \end{tabular}
    \label{tab:mc_l}
\end{table}

\begin{table}[!ht]
    \centering
    \small
    \caption{Numerical Results on More Large-Scale Matrix Completion Problems}
    \vspace{5pt}
    \begin{tabular}{|cc|cc||cc|cc|}
    \hline
        \multicolumn{2}{|c|}{\textbf{~}} & \multicolumn{2}{c||}{\textbf{cuLoRADS}} & \multicolumn{2}{c|}{\textbf{~}} & \multicolumn{2}{c|}{\textbf{cuLoRADS}} \\ 
        n & m & time & error & n & m & time & error \\ \hline
        60000 & 23683563 & 2.79 & 5.51E-07 & 200000 & 7996791 & 0.44 & 1.39E-06 \\ 
        80000 & 31682246 & 5.56 & 1.64E-06 & 300000 & 11996886 & 0.65 & 1.13E-05 \\ 
        100000 & 39682090 & 7.43 & 2.19E-06 & 400000 & 15996836 & 1.84 & 3.61E-06 \\ 
        120000 & 47682387 & 7.4 & 4.17E-07 & 600000 & 23996846 & 1.74 & 1.46E-07 \\ 
        140000 & 55682003 & 8.24 & 4.25E-07 & 800000 & 31996932 & 2.52 & 8.66E-06 \\ 
        160000 & 63681433 & 10.09 & 2.58E-07 & 1000000 & 39996783 & 2.61 & 2.38E-07 \\ 
        180000 & 71681424 & 15.67 & 4.53E-06 & 1200000 & 47996814 & 3.31 & 5.07E-06 \\ 
        200000 & 79681236 & 7.38 & 1.59E-06 & 1400000 & 55996812 & 4.23 & 2.12E-06 \\ 
        300000 & 119680181 & 8.35 & 9.79E-07 & 1600000 & 63996921 & 6.21 & 1.68E-06 \\ 
        400000 & 159679797 & 11.38 & 1.81E-06 & 1800000 & 71996873 & 11.2 & 5.93E-09 \\ 
        600000 & 239681020 & 19.25 & 5.63E-07 & 2000000 & 79996838 & 5.93 & 4.19E-07 \\ 
        800000 & 319680602 & 39.69 & 5.68E-07 & 4000000 & 159996869 & 63.45 & 2.09E-06 \\ 
        ~ & ~ & ~ & ~ & 6000000 & 239996782 & 88.67 & 3.11E-08 \\ 
        ~ & ~ & ~ & ~ & 8000000 & 319996871 & 85.83 & 4.98E-07 \\ 
        ~ & ~ & ~ & ~ & 20000000 & 199999806 & 442.55 & 8.05E-07 \\ \hline
    \end{tabular}
    \label{tab:mc_h}
\end{table}

\subsection{Capability on Solving General SDP Benchmark Problems\label{sec:benchmark_results}}
In this section, we benchmark the performance of cuLoRADS on a set of general semidefinite programming problems from Hans Mittelmann's SDP benchmarks\footnote{\url{https://plato.asu.edu/ftp/sparse_sdp.html}}. The experiment showcases the GPU solver's performance across various types of SDP problems, though their scales are relatively small. Since cuLoRADS is still in its early stages of development and does not yet support solving problems with linear variables, we include only a subset of the benchmark problems that do not have this feature.

Table~\ref{tab:general} demonstrates the performance of cuLoRADS on several different types of SDP problems compared to other benchmark solvers. cuLoRADS exhibits significantly faster solving times than LoRADS by leveraging the computational power of GPUs and the acceleration design. The quality of the solutions obtained is also superior to that of LoRADS due to the re-opt technique.

Moreover, cuLoRADS performs competitively against three other state-of-the-art solvers, solving the majority of problems faster than the others on GPU. Notably, cuLoRADS also shows strong performance compared to the commercial solver COPT on some SDP problems, such as AlH, where it is 200 times faster. It consistently achieves shorter solving times on other problems listed in Table~\ref{tab:general}. These results highlight the substantial potential of GPUs in accelerating SDP solving, even when considering commercial-grade performance.
\begin{table}[!ht]
    \centering
    \scriptsize
    \caption{Numerical Results on General SDP Problems}
    \vspace{5pt}
    \setlength{\tabcolsep}{1.2mm}
    \begin{tabular}{|ccc|cc|cc|cc|cc|cc|}
    \hline
        ~ & ~ & ~ & \multicolumn{2}{c|}{\textbf{cuLoRADS}} & \multicolumn{2}{c|}{\textbf{LoRADS}} & \multicolumn{2}{c|}{\textbf{SDPLR}} & \multicolumn{2}{c|}{\textbf{SDPNAL+}} & \multicolumn{2}{c|}{\textbf{COPT}} \\ 
        name & n & m & time & $\text{error}_{\text{max}}$ & time & $\text{error}_{\text{max}}$ & time & $\text{error}_{\text{max}}$ & time &  $\text{error}_{\text{max}}$ & time &  $\text{error}_{\text{max}}$ \\ \hline
        1dc.1024 & 1024  & 24064  & 47.3  & 8.96E-06 & 792  & 1.69E-04 & 6275  & 6.23E-05 & f & f & 232 & 2.04E-08 \\ 
        1et.2048 & 2048  & 22529  & 68.5  & 1.42E-05 & 1811  & 9.50E-05 & f & f & 1481  & 2.30E-05 & 261 & 5.36E-08 \\ 
        1tc.2048 & 2048  & 18945  & 71.3  & 1.43E-05 & 1660  & 5.05E-05 & f & f & 3724  & 4.80E-06 & 198 & 4.75E-08 \\ 
        1zc.1024 & 1024  & 16641  & 20.3  & 9.93E-06 & 608  & 2.35E-04 & 914  & 2.54E-05 & 57.6  & 6.90E-07 & 72.4  & 2.98E-08 \\ 
        AlH & 5990  & 7230  & 8.8  & 1.64E-05 & 63.2 & 1.23E-03 & 2491 & 4.87E-04 & 943 & 6.70E-06 & 1208 & 2.51E-09 \\ 
        BH2 & 2166  & 1743  & 7.3  & 2.66E-05 & 7.8 & 1.82E-02 & 124 & 1.24E-04 & 48.3 & 1.78E-06 & 35.7 & 2.54E-09 \\ 
        cancer\_100 & 569  & 10470  & 35.6  & 8.63E-06 & 514 & 3.28E-04 & t & t & 728 & 4.43E-05 & 83.4 & 1.97E-05 \\ 
        CH2 & 2166  & 1743  & 9.7  & 1.79E-05 & 6.7 & 2.09E-03 & 51.5 & 2.49E-04 & 88.8 & 3.48E-04 & 19.1 & 3.56E-09 \\ 
        checker\_1.5 & 3970  & 3971  & 1.2  & 1.61E-07 & 1.0 & 1.54E-07 & 40.5 & 1.36E-04 & 6201 & 1.25E-03 & 20.3 & 4.53E-08 \\ 
        cphil12 & 364  & 12376  & 49.7  & 1.91E-09 & f & f & 0.4  & 3.33E-03 & 0.2  & 2.68E-16 & 15.0  & 1.58E-11 \\ 
        fap09 & 174  & 30276  & 56.7  & 3.89E-05 & f & f & 574  & 9.85E-05 & f & f & 53.2  & 2.99E-07 \\ 
        foot & 2208  & 2209  & 10.7  & 5.52E-06 & 1198  & 1.21E-02 & f & f & 2777  & 2.09E-05 & 8.6  & 4.92E-06 \\ 
        G40\_mb & 2000  & 2001  & 3.9  & 2.26E-07 & 19.9 & 1.29E-07 & 60.0 & 1.07E-06 & 4645 & 3.33E-07 & 12.0 & 1.07E-05 \\ 
        G48\_mb & 3000  & 3001  & 2.6  & 3.26E-06 & 17.2 & 1.37E-07 & 25.5 & 2.96E-04 & 3019 & 7.15E-02 & 17.5 & 4.85E-03 \\ 
        G60\_mb & 7000  & 7001  & 27.9 & 2.47E-06 & 2271  & 1.14E-05 & 985  & 8.51E-07 & t & t & 105 & 8.83E-08 \\ 
        H3O & 3162  & 2964  & 4.6  & 2.05E-05 & 10.1 & 2.47E-04 & 199 & 7.14E-05 & 295 & 5.39E-06 & 101 & 3.49E-09 \\ 
        hand & 1296  & 1297  & 2.1  & 4.02E-07 & 116 & 7.62E-08 & 24.2 & 1.08E-07 & 628 & 2.14E-07 & 92.8 & 4.09E-05 \\ 
        ice\_2.0 & 8113  & 8113  & 3.5  & 4.22E-06 & 4.9 & 2.60E-06 & 42.0 & 1.01E-05 & t & t & 269 & 1.98E-08 \\ 
        NH2 & 2046  & 1743  & 2.5  & 2.41E-05 & 24.8 & 1.68E-02 & 19.0  & 8.48E-05 & 10.8  & 9.10E-06 & 18.4 & 2.51E-09 \\ 
        NH3 & 3162  & 2964  & 3.2  & 3.53E-05 & 9.6 & 1.60E-04 & 177 & 2.42E-04 & 55.4 & 2.65E-05 & 63.8 & 2.98E-09 \\ 
        NH4 & 4426  & 4743  & 2.1  & 2.49E-05 & 27.1  & 8.95E-02 & 231  & 3.12E-04 & 105  & 7.47E-04 & 329 & 4.97E-09 \\ 
        p\_auss2\_3.0 & 9115  & 9115  & 4.1  & 5.40E-06 & 4.4 & 9.42E-07 & 84.0 & 3.14E-05 & t & t & 289 & 4.74E-06 \\ 
        rendl1 & 2000  & 2001  & 9.4  & 2.74E-06 & 1576  & 6.29E-05 & 306  & 1.53E-08 & 638  & 4.12E-07 & 4.7  & 1.38E-05 \\ 
        sensor\_1000 & 6528  & 5549  & 37.6  & 6.86E-06 & f & f & 719  & 9.87E-03 & 739  & 2.53E-02 & 66.6  & 1.03E-05 \\ 
        theta12 & 600  & 17979  & 18.6  & 5.85E-06 & 86.2 & 1.80E-03 & 127 & 1.83E-05 & 11.7 & 5.31E-06 & 9.9 & 1.74E-08 \\ 
        theta102 & 500  & 37467  & 19.4  & 2.95E-05 & 71.4 & 3.02E-03 & 141 & 2.53E-05 & 6.4 & 1.12E-06 & 4.6 & 2.80E-09 \\ 
        theta123 & 600  & 90020  & 25.9  & 8.16E-06 & 206 & 8.26E-04 & 369 & 3.50E-05 & 10.6 & 3.86E-07 & 8.9 & 1.90E-08 \\ \hline

    \end{tabular}
    \label{tab:general}
\end{table}

\section{Conclusion\label{sec:conclusion}}
In this paper, we tackle the long-standing challenge of developing a scalable and efficient solver for semidefinite programming problems. We propose accelerating low-rank factorization-based methods using GPUs, providing solutions to the efficiency and scalability bottlenecks on GPUs. We also introduce cuLoRADS, a GPU-based solver implementing the low-rank splitting ADMM approach described in \cite{han2024low}. By incorporating both the classic Burer-Monteiro method and ADMM, cuLoRADS demonstrates extraordinary performance across a range of SDP problems. Our results highlight its remarkable efficiency and scalability, significantly outperforming traditional CPU solvers.

\paragraph{Limitations and future works.}

In its current stage of development, the LoRADS approach may still exhibit inefficiencies when applied to certain types of problems. Efforts are ongoing to enhance its robustness across a broader spectrum of general semidefinite problems. Notably, our acceleration techniques are compatible with many low-rank factorization-based methods, provided their computations involve basic matrix-matrix, matrix-vector, and vector-vector operations. For example, algorithms such as the one proposed by \cite{monteiro2024low} can be explored for potential acceleration using our techniques. Additionally, we are focusing on kernel-level optimization for the custom kernels employed in our calculations, as well as investigating the potential of multi-GPU implementations and mixed-precision computation.

\bibliographystyle{plainnat}
\bibliography{ref}

\newpage
\appendix

\section{Detailed Algorithm Flow of LRALM}

\subsection{Dynamic Rank Updating. }
Following the approach described by~\cite{han2024low}, we use $r = O(\log m)$ as the initial rank, which is recommended by~\cite{so2007theory}. Moreover, we adopt a progressive rank update method to achieve better algorithm performance. 
The timing for increasing the rank is determined heuristically: if the number of iterations of LBFGS (Algorithm~\ref{alg:lbfgs} for solving the subproblem) exceeds a predetermined threshold, we adjust the rank. Specifically, the rank is multiplied by a factor of 1.5 (i.e., $\lceil 1.5 r\rceil$), with the additional elements initialized to values sampled from a normal distribution or set to zero during the rank update.
Additionally, since the maximum rank is capped at $\sqrt{2m}$, we take the minimum of the two values after updating the rank. Hence, the rank updating can be summarized as $r_{\text{new}} = \min\bcbra{\lceil 1.5 r_{\text{old}}\rceil, \sqrt{2m}}$. In practice, we use a conservative threshold for rank updating and can typically solve the problem at the $O(\log m)$ rank order.

\subsection{Exact Line Search. \label{sec:exact_line_search}}

For clearness, we consider how to conduct an exact line search for a single block variable. Similar conclusions can be directly extended to the multi-block case. The following deduction is completely consistent with paper \cite{BMImplement}. We rewrite it here for completeness. The exact line search aims to find the optimal step size $\tau$ that minimizes the following function:
\begin{equation}
\begin{split}\mcal L(R+\tau D)= & \inner C{\rbra{R+\tau D}\rbra{R+\tau D}^{\top}}+\lambda^{\top}\Brbra{b-\mcal A\brbra{\rbra{R+\tau D}\rbra{R+\tau D}^{\top}}}\\
 & +\frac{\rho}{2}\Bnorm{b-\mcal A\brbra{\rbra{R+\tau D}\rbra{R+\tau D}^{\top}}}^{2}\\
= & \inner C{RR^{\top}}+\tau\inner C{RD^{\top}+DR^{\top}}+\tau^{2}\inner C{DD^{\top}}\\
 & +\lambda^{\top}\Brbra{b-\mcal A\brbra{RR^{\top}+\tau^{2}DD^{\top}+\tau RD^{\top}+\tau DR^{\top}}}\\
 & +\frac{\rho}{2}\Bnorm{b-\mcal A\brbra{RR^{\top}+\tau^{2}DD^{\top}+\tau RD^{\top}+\tau DR^{\top}}}^{2}.
\end{split}
\end{equation}
Function $\mcal L(R+\tau D)$ is a polynomial of $\tau$ and discard
some constant, then we have the function 
\begin{equation*}
\begin{aligned}
\mcal L(R+\tau D)\sim & \tau\inner C{RD^{\top}+DR^{\top}}+\tau^{2}\inner C{DD^{\top}}\\
&+\lambda^{\top}\Brbra{b-\mcal A\brbra{RR^{\top}+\tau RD^{\top}+\tau DR^{\top}+\tau^{2}DD^{\top}}}\\
 & +\frac{\rho}{2}\Bnorm{b-\mcal A\Brbra{RR^{\top}+\tau RD^{\top}+\tau DR^{\top}+\tau^{2}DD^{\top}}}^{2}\\
\sim & \tau\inner C{RD^{\top}+DR^{\top}}+\tau^{2}\inner C{DD^{\top}}+\lambda^{\top}\Brbra{b-\mcal A\brbra{RR^{\top}+\tau RD^{\top}+\tau DR^{\top}+\tau^{2}DD^{\top}}}\\
 & +\rho\brbra{b-\mcal A(RR^{\top})}^{\top}\mcal A\Brbra{\tau RD^{\top}+\tau DR^{\top}+\tau^{2}DD^{\top}}\\
 & +\frac{\rho\tau^{2}}{2}\Brbra{\norm{\mcal A\brbra{RD^{\top}+DR^{\top}}}^{2}+\tau^{2}\norm{\mcal A(DD^{\top})}^{2}+2\tau\brbra{\mcal A\brbra{RD^{\top}+DR^{\top}}^{\top}\mcal A(DD^{\top})}}\\
= & \tau\Brbra{\inner C{RD^{\top}+DR^{\top}}-\brbra{\lambda+\rho\brbra{b-\mcal A(RR^{\top})}}^{\top}\mcal A(RD^{\top}+DR^{\top})}\\
 & +\tau^{2}\Brbra{\inner C{DD^{\top}}-\brbra{\lambda+\rho\brbra{b-\mcal A(RR^{\top})}^{\top}\mcal A(DD^{\top})}+\frac{\rho}{2}\brbra{\norm{\mcal A(RD^{\top}+DR^{\top})}^{2}}}\\
 & +\tau^{3}\Brbra{\rho\mcal A(RD^{\top}+DR^{\top})^{\top}\mcal A(DD^{\top})}\\
 & +\tau^{4}\cdot\frac{\rho}{2}\norm{\mcal A(DD^{\top})}^{2},
\end{aligned}
\end{equation*}
where symbol $"\sim"$ means that two functions have the same minimizers. 
Hence, the optimal $\tau^{*}$ can be obtain by solving the following
equation:
\begin{equation*}
4a_1\tau^{3}+3a_2\tau^{2}+2a_3\tau+a_4=0,
\end{equation*}
where 
\begin{equation*}
a_1:=\frac{\rho}{2}\norm{q_{2}}^{2},a_2=\rho q_{1}^{\top}q_{2},a_3=p_{2}-(\lambda+\rho q_{0})^{\top}q_{2}+\frac{\rho}{2}\norm{q_{1}}^{2},a_4=p_{1}-(\lambda+\rho q_{0})^{\top}q_{1}
\end{equation*}
and 
\begin{align*}
p_{1} & =\inner C{RD^{\top}+DR^{\top}},p_{2}=\inner C{DD^{\top}},q_{0}=b-\mcal A(RR^{\top})\\
q_{1} & =\mcal A(RD^{\top}+DR^{\top}),q_{2}=\mcal A(DD^{\top}).
\end{align*}

\end{document}